 \documentclass[11pt, a4paper]{article}

 \pdfoutput=1

 \usepackage{amsmath}
 \usepackage{amssymb}
 \usepackage{amsthm}
 \usepackage{amsfonts}
 \usepackage{graphicx}
 \usepackage{enumerate}
 \usepackage{mathrsfs}
 \usepackage{bbm}

 \usepackage{hyperref}
 \usepackage{fullpage}
 \usepackage{xypic}
 
 \usepackage{color}
 
 \usepackage[nottoc,numbib]{tocbibind}
 
 \usepackage{cleveref}
 
 \usepackage[perpage]{footmisc}

 \newtheorem{thm}{Theorem}[section]
 
 \newtheorem{lem}[thm]{Lemma}

 \theoremstyle{definition}
 \newtheorem{exm}[thm]{Example}
 
 \newtheorem{rmk}[thm]{Remark}

 \newtheorem{claim}[thm]{Claim}  

 \newtheoremstyle{nonum}{}{}{\itshape}{}{\bfseries}{.}{ }{\thmnote{#3}}
 \theoremstyle{nonum}

 \DeclareMathOperator{\supp}{supp}

 \DeclareMathOperator{\codim}{codim}

  \DeclareMathOperator{\Vol}{Vol}

  \DeclareMathOperator{\Area}{Area}
  \DeclareMathOperator{\diam}{diam}

 \newcommand{\tit}{\textit}

 \newcommand{\mrm}{\mathrm}

 \newcommand{\calF}{{\mathcal{F}}}

 \newcommand{\eps}{\varepsilon}
 \newcommand{\epsi}{\varepsilon}

 \newcommand{\Sig}{\Sigma}

 \newcommand{\ph}{\varphi}
 \newcommand{\vphi}{\varphi}
 \newcommand{\om}{\omega}
 
 \newcommand{\na}{\nabla}

 \newcommand{\R}{\mathbb R}

 \newcommand{\N}{\mathbb N}

 \newcommand{\setm}{\setminus}

 \newcommand{\pa}{\partial}
 \newcommand{\de}{\partial}

 \newcommand{\restr}{\big|}  
 
 \newcommand{\Czto}{\displaystyle \xrightarrow[C_0]{}}

 \newcommand{\pois}[1]{\{#1\}} 
 \DeclareMathOperator{\sgrad}{sgrad}
 
 \newcommand{\Symp}{\mrm{Symp}}

 \newcommand{\poisnm}[1]{\|\{#1\}\|} 

 \newcommand{\ol}{\overline}
 \newcommand{\til}{\tilde}

 \newcommand{\into}{\hookrightarrow}

 \newcommand{\vect}[1]{\overrightarrow{#1}}

\begin {document}

\title{Rigidity versus flexibility of the Poisson bracket \\ with respect to the $L_p$-norm}
\author{Karina Samvelyan 
		\footnote{Partially supported by the European Research Council Advanced grant 338809.}
		}
	
\maketitle


\begin{abstract}
	Rigidity of the Poisson bracket with respect to the uniform norm is one of the central phenomena discovered within function theory on symplectic manifolds. In the present work we examine the case of $L_p$ norms with
	$p < \infty$. We show that $L_p $ - Poisson bracket invariants exhibit rigid behavior in dimension two, and we provide an evidence for their flexibility in higher dimensions.
	
\end{abstract}

\pagenumbering{arabic}

\tableofcontents 


\section{Introduction and statement of results}

The subject of the present work is function theory on symplectic manifolds.
We focus on the interplay between rigidity and flexibility of the Poisson bracket.\\

Recall that a symplectic structure on an even-dimensional manifold $M^{2n}$ is a closed differential 2-form $\om$, whose top power $\om^n$ vanishes nowhere.
The classical Darboux theorem states that locally any symplectic manifold looks as the standard symplectic vector space $\R^{2n}$ with coordinates $(p_1, \ldots, p_n, q_1, \ldots, q_n)$ equipped with the symplectic form $\sum_{i=1}^n dp_i \wedge dq_i$.
Another important example of a symplectic manifold is a surface equipped with an area form.

A fundamental notion of symplectic geometry is the Poisson bracket, $\pois{F,G}$, of a pair of smooth functions $F$ and $G$ on $M$. Locally, in Darboux coordinates $p_i, q_i$ ($i = 1\ldots n$),
\begin{align}
\pois{F,G} = \sum_{i=1}^n
\left( 
\frac{\de F}{\de q_i}\frac{\de G}{\de p_i} - \frac{\de G}{\de q_i}\frac{\de F}{\de p_i}
\right) \;.
\end{align}

\noindent
The following identity can provide a coordinate-free definition:
\begin{align} \label{poisson_brck_def_1}
\pois{F,G} \om^n = -n \cdot dF \wedge dG \wedge \om^{n-1} \;.
\end{align}

%

\subsection{Measurements with the Poisson bracket}

Let $(M^{2n}, \om)$ be a symplectic manifold. A significant character of our story is the functional
\begin{align*}
\Phi_p : C^\infty_c(M) \times C^\infty_c(M) \to \R_{\geq 0},\ 
(F,G) \mapsto \| \pois{F,G} \|_p \;.
\end{align*}

\noindent
Here $C_c^\infty (M)$ stands for the space of smooth compactly supported functions on $M$, and we write $\| F\|_p$ for the $L_p$-norm 
$$ 
\| F \|_p = \left( \int_M |F|^p \om^n \right)^{1/p}
$$
associated to the volume form $\om^n$ on $M$.
We consider $p \in [1,\infty]$, where by $L_\infty$-norm we mean the uniform norm 
$\| F \|_\infty = \max_M |F|$.\\

It was shown 
that for $p=\infty$ this functional, $\Phi_\infty$, is lower semi-continuous with respect to the $L_\infty$-norm on $C^\infty_c (M)$. (See \cite{polterov_rosen2014func_theory}, \cite{entov_polterov_c0_rig_poiss_br} and \cite{buhovsky_2/3_conv_pb_2010}. These texts deal with the multidimensional case, extending previous results by Cardin-Viterbo (\cite{cardin_viterbo_08}) and Zapolsky (\cite{zapolsky_quasi_and_pb_surfaces_07}).)
This fact is quite surprising, since the Poisson bracket depends on the first derivatives of the functions, while the convergence is in the uniform norm only. 
Let us mention also, that the functional $\Phi_p$ is not continuous, as we can slightly alter the two functions, while changing their derivatives extensively. \\

\noindent
Our first result deals with the behaviour of the functional $\Phi_p$ in the $L_q$-topology for general  $p$ and $q$.

\begin{thm} \label{thm: thm_pb_non_rigidity_L_p_L_q_combined_intro}
	Let $1 \leq q <\infty$ and $1\leq p \leq \infty$. For any two functions $F,G \in C^\infty_c(M)$ that are not Poisson commuting ($\pois{F,G} \neq 0$),
	there exist two sequences $F_N,\,G_N \in C^\infty_c (M)$ with $F_N \xrightarrow[{L_q}]{} F$, $G_N \xrightarrow[{L_q}]{} G$ and $\poisnm{F_N,G_N} _p \to 0$ as $N \to \infty$. 
	
\end{thm}


\noindent
In fact, we will construct two sequences satisfying
$F_N \xrightarrow[{L_\infty}]{} F$, $G_N \xrightarrow[{L_q}]{} G$ 
with ${\pois{F_N,G_N} \equiv 0}$.\\

Thus, in these cases the semicontinuity phenomenon disappears and the rigidity we witnessed in the case of the uniform norm is replaced by flexibility.\\
The case $q = \infty$, $p < \infty$ remains open.

\subsection{Poisson bracket invariant of quadruples, $pb_4^q$}

Next, we discuss another measurement that has to do with the Poisson bracket.
Let $X_0, X_1, Y_0, Y_1$ be compact subsets of a symplectic manifold $(M,\om)$, such that 
$X_0 \cap X_1 = Y_0 \cap Y_1 = \emptyset$.
Fix some $1\leq q \leq \infty$ and set

\begin{equation}\label{eq:intro_pb_4_def_1}
pb_4 ^q (X_0,X_1,Y_0,Y_1) = \inf_{(F,G)} \| \pois{F,G} \|_q \;,
\end{equation}
where the infimum is taken over all pairs $F,G\in C^\infty_c(M)$,
such that 
$$
F \restr_{X_0} \leq 0,\ F \restr_{X_1} \geq 1 \;,\ 
G \restr_{Y_0} \leq 0,\ G \restr_{Y_1} \geq 1 \;.
$$
In the notation $pb_4^q$, $pb$ stands for Poisson bracket, the subindex 4 is for the fact that we deal with a quadruple of subsets, and $q$ is to signify the $L_q$-norm.
\begin{figure} [h!]
	\centering
	\includegraphics[scale = 0.75, trim = 0 0 0 0 ]{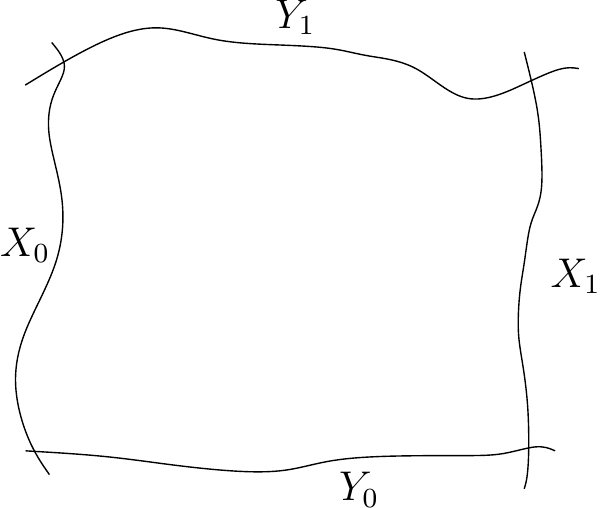}
\end{figure}

It is known (see \cite{buhovsky_entov_polterovich_2012poisson_brck}, \cite{polterov_rosen2014func_theory}) that for certain quadruples of subsets, $pb_4^\infty$ is strictly positive, thus manifesting the rigidity of the uniform norm of the Poisson bracket.
In this work, we explore the properties of the functional $pb_4^q$ also when $1 \leq q < \infty$.

We show that if $\dim M = 2$, i.e. in the case of $M$ being a surface, rigidity of $pb_4^q$ persists, whereas in the multidimensional case $pb_4^q$ exhibits flexible behavior.

\subsubsection{Rigidity in the 2-dimensional case}

We shall consider the invariant $pb_4^q$ of the four sides of a quadrilateral 
on a smooth surface $M$ equipped with an area form $\om$.
%
For us, a curvilinear quadrilateral on a smooth surface $M$ is the image of an embedding of a square $[0,1]^2 \subset \R^2$ into the interior of $M$. 

Suppose that $X_0,Y_0, X_1,Y_1$ are sides of a curvilinear quadrilateral $\Pi \subset M$ taken in counterclockwise order.
We consider $pb_4^q (\Pi) := pb_4^q (X_0,X_1,Y_0,Y_1)$.

\noindent
It turns out that in the case $q>1$ the value of $pb_4^q(\Pi)$ depends on the areas of $\Pi$ and $M$, while $pb_4^{q=1} (\Pi)$ is independent of these areas.
%

%
%
%
%


\begin{thm} \label{thm: intro_thm_pb_4_dim=2_all}
	Let $(M, \om)$ be a connected symplectic surface without boundary and let $\Pi \subset M$ be a curvilinear quadrilateral.\\
	Let $1 \leq q < \infty$. Denote $A = \Area_\om(\Pi)$, $B = \Area_\om(M)$.\\ 
	(i) If $\Area(M) < \infty$, then 
	$\displaystyle{pb_4^q (\Pi) = \Big( \frac{1}{A^{q-1}} + \frac{1}{(B-A)^{q-1}} \Big)^{1/q}}\ $.
	In particular, $pb_4^1 (\Pi) = 2$.\\
	(ii) If $\Area(M) = \infty$, then 
	$\displaystyle{pb_4^q(\Pi) =  \Big( \frac{1}{A^{q-1}} \Big)^{1/q}}\ $. 
\end{thm}

\begin{rmk}
	Note that (ii) is a limiting case of (i) as $B \to \infty$.
\end{rmk}

\subsubsection{$pb_4^q$ of a curve on a surface}

The quantity $pb_4^q$ gives rise to an invariant of simple closed curves on surfaces. Consider such a curve $\tau$ on a smooth connected oriented surface $\Sigma$ without boundary. \\
Divide the curve into four segments $\Delta_1, \Delta_2, \Delta_3, \Delta_4$ and consider $pb_4^q (\Delta_1, \Delta_2, \Delta_3, \Delta_4)$ of this quadruple. 
We will see that this quantity does not depend on the division of $\tau$, and thus this construction
defines an invariant $pb_4^q(\tau)$ of the curve $\tau$. To the best of our knowledge, this definition is new even for $q=+\infty$.\\
It appears that $pb_4^q (\tau)$ captures some topological information regarding the curve $\tau$. 
Namely, it distinguishes separating simple closed curves from non-separating ones. 
Recall that $\tau$ is called non-separating if $\Sigma \setm \tau$ is connected.
If a curve is non-separating, $pb_4^q(\tau)$ vanishes, while it is not the case for a separating curve.

\begin{thm}	\label{thm: intro_thm_pb_4_non_sep_curve_on_surface}
	Let $(\Sigma,\om)$ be a smooth connected symplectic surface without boundary, 
	and let $\tau \subset \Sigma$ be a smooth simple closed curve.
	If $\tau$ is non-separating, then $pb_4^q (\tau) = 0$ for any $1 \leq q \leq \infty$.
\end{thm}


\begin{thm} \label{thm: intro_pb_4_sep_curve_surface_no_bdry_finite_areaed_components}
	Let $(\Sigma,\om)$ be a smooth connected symplectic surface without boundary, 
	and let $\tau \subset \Sigma$ be a smooth simple closed separating curve. 
	Suppose that the components $\Sig_1$ and $\Sig_2$ of $\Sigma \setm \tau$ have finite areas $A$ and $B$ respectively.
	Then $pb_4^q$ does not vanish, and moreover,
	\begin{equation}
	pb_4^q(\tau) = 	
	\begin{cases}
	2							 & \text{if } \ q=1  \;, \\
	\left(\frac{1}{A^{q-1}} + \frac{1}{B^{q-1}} \right) ^{1/q}       & \text{if } \ 1< q < \infty \;, \\
	\max (\frac{1}{A}, \frac{1}{B})				     & \text{if } \ q=\infty \;.
	\end{cases}	
	\end{equation}
\end{thm}

\subsubsection{$pb_4^q$: the multidimensional case}


Here we present a new mechanism revealing that $pb_4^q$ vanishes in higher dimensions in certain situations.


\begin{thm} \label{thm: intro_thm_pb_4_vanish_high_dim}
	Let $X_0, X_1, Y_0, Y_1$ be compact subsets of a symplectic manifold $(M^{2n},\om)$, where $X_1$ is a submanifold, such that 
	$X_0 \cap X_1 = Y_0 \cap Y_1 = \emptyset$. 
	Denote $d = \dim X_1$ and suppose $d \leq 2n-2$. Then 
	$pb_4^q(X_0, X_1, Y_0, Y_1) = 0$ whenever $q \leq 2n - d, n\geq 2$.
\end{thm}

Interestingly enough, $pb_4^\infty$ for such a quadruple can be positive. For instance, examine
$$
{[0,1]^2 \times T^* S^1 \subset  \R^2 \times T^* S^1}
$$
and denote the sides of $[0,1]^2$
by $a, b, c, d$, listed in cyclic order. Pick a fixed circle (the zero section) $S^1$ on the cylinder $T^* S^1$.
Consider the quadruple 
$$
(X_0, Y_0, X_1, Y_1) = (a \times S^1, b \times S^1, c \times S^1, d \times S^1) \;,
$$
which is called \emph{the stabilization} of $(a,b,c,d)$ 
(see \cite{buhovsky_entov_polterovich_2012poisson_brck}). 
Here for $X_1 = c \times S^1$, we have $d = \dim X_1 = 2$ and $n=2$. 
For $q = \infty$ we have 
$pb_4^\infty (a \times S^1, b \times S^1, c \times S^1, d \times S^1) > 0$,
see \cite[section 7.5.4]{polterov_rosen2014func_theory},
i.e. positivity of $pb_4^\infty$ on the sides of the quadrilateral survives the stabilization. 
\cref{thm: intro_thm_pb_4_vanish_high_dim} above shows that this is not longer valid for $q \leq 2n-d=2$. 
The case of finite $q > 2$ is currently out of reach.







\section{Poisson Bracket flexibility with respect to $L_p$-norms}


Let $(M^{2n},\om)$ be a symplectic manifold, $n\in \N$, and fix $1\leq p \leq \infty$ and $1 \leq q < \infty$. 
Denote by $C^\infty_c(M)$ the space of smooth functions on $M$ with compact support.

\begin{thm} \label{thm: thm_pb_non_rigidity_L_p_L_q_combined}
	%
	For any two functions $F,G \in C_c^\infty (M)$ there exist two sequences 
	${F_N,\,G_N \in C_c^\infty (M)}$ with $F_N \xrightarrow[{C^0}]{} F$, $G_N \xrightarrow[{L_q}]{} G$ and 
	$\pois{F_N,G_N} = 0\ $  $\forall N \in \N$.
\end{thm}

\begin {proof}
Let us note first that, in the notations of \cref{thm: thm_pb_non_rigidity_L_p_L_q_combined}, 
obtaining 
$F_N \Czto F$ 
would be sufficient to deduce \cref{thm: thm_pb_non_rigidity_L_p_L_q_combined_intro} for any $1 \leq p \leq \infty$, as long as all the functions $F_N$ will be supported on a compact set independent of $N$, which indeed will be the case in our construction below.\\
Given any non-commuting $F,G \in C^\infty_c (M)$, we shall construct $\til{F}$ and $\til{G}$ with $\pois{\til{F},\til{G}} = 0$, such that they are arbitrarily close to $F$ and $G$ in the norms $C_0$ and $L_p$ respectively.\\

Let us fix some Riemannian metric $d$ on $M$. We will only deal with a compact subset of $M$ (where our functions will be supported), any two metrics on this compact are equivalent, so the choice of metric would not effect our argument.\\




By a \tit{simplex in $M^{2n}$} we mean the image of an embedding $\Delta \to M$, where $\Delta$ is a (closed) simplex in $\R ^{2n}$. 
A triangulation of $M$ is a representation of $M$ as a union of such simplices. We also require each two simplices to intersect only in a common face, which is a simplex of lower dimension.\\
A construction described in \cite{cairns1961} produces such a triangulation of $M$, representing it as a finite union a simplices (for a compact $M$).
In case of a non-compact manifold, we will only need a triangulation of $\supp(F) \cup \supp (G)$.
Moreover, using the same procedure, we can make the diameter of all simplices to be smaller than any prescribed constant. (Here the diameter is with respect to the chosen metric $d$.) \\

Let $\eps > 0$. Take such a triangulation (of $\supp(F) \cup \supp(G) $) with all simplices having diameter $< \delta$, where $\delta > 0$ will be fixed later and will depend on $\eps$ and $F$.
%
Note that given a simplex $Q$ from this triangulation, we can find an open subset $Q' \subseteq Q$, such that $Q \setm Q' \supseteq \pa Q$ and $\Vol (Q \setm Q') \leq a \cdot \Vol (Q)$ for a (small) fixed $a>0$ (i.e. $Q'$ occupies most of the volume of $Q$.) By $\Vol$ here and later in the proof we mean volume with respect to $\om^n$. 


For every simplex $Q$ from the triangulation of $M$, we shall take open subsets with smooth boundary 
$Q_3 \Subset Q_2 \Subset Q_1 \Subset Q$, satisfying 
$\Vol (Q \setm Q_3) \leq \Vol(Q) \cdot \eps$ (here $A \Subset B$ means $Cl(A) \subseteq int(B)$).
This last condition will be essential for taking a suitable $\til{G}$.
%
%

\subsubsection*{Construction of $\til{F}$}

Consider a simplex $Q$ with open subsets $Q_2 \Subset Q_1 \Subset Q$. We take an auxiliary smooth function $\ph : Q \to [0,1]$ such that $\ph \restr_{Q_2} \equiv 0$ and $\ph \restr_{Q \setm Q_1} \equiv 1$. Fix also a point $x_0 \in Q_2$. Define $\til{F}$ on $Q$ to be

$$\til{F} (x) = \ph (x) F(x) + (1 - \ph (x)) F(x_0) \;.$$

\noindent
So on $Q_2$ we have $\til{F} \equiv F(x_0)$ ($\til{F}$ being an approximation of $F$ on $Q_2$), while outside $Q_1$, $\til{F} \equiv F$. (See \cref{fig: construction_F_G_tildes_PB_flexibility}.)

Next, glue all these $\til{F}$ hereby defined on each simplex. It is possible, since on adjacent simplices, in a neighborhood of their intersection the patches of $\til{F}$ are equal to $F$.
We get a compactly supported smooth function $\tilde{F}$ on $M$, as $F$ is compactly supported. 
%
%
%
Note also that $F$ is uniformly continuous on its (compact) support, i.e. for any $\epsi > 0$ there exists some $\delta > 0$, so that $ d(x,y) < \delta$ (Riemannian distance) implies $| F(x) - F(y) | < \epsi$.

\noindent
Thus, taking appropriate $\delta > 0$, on a single simplex $Q$, for any $x \in Q $ we have 
\begin{align*}
|\til{F} (x) - F(x)| & = |\ph(x)F(x) + (1-\ph(x))F(x_0) - F(x)| =\\
& =  \underbrace{|1-\ph(x)|}_{\leq 1} \cdot \underbrace{|F(x)-F(x_0)|}_{< \epsi} 
< \eps \;,
\end{align*}
where the last inequality hold since $\diam (Q) < \delta$.
So $\| F-\til{F} \| _\infty \leq \eps$ on each $Q$ 
taking $\delta > 0$ small enough to suite all simplices. 
Hence $\| F-\til{F} \| _\infty \leq \eps$ on the whole $M$.
Thus, $\| F-\til{F} \| _\infty$ 
and, consequently, $\| F-\til{F} \| _q$ can be made as small as we wish, taking $\delta > 0$ small enough.

\begin{figure} [h!]
	\centering
	\includegraphics[scale= 0.9]{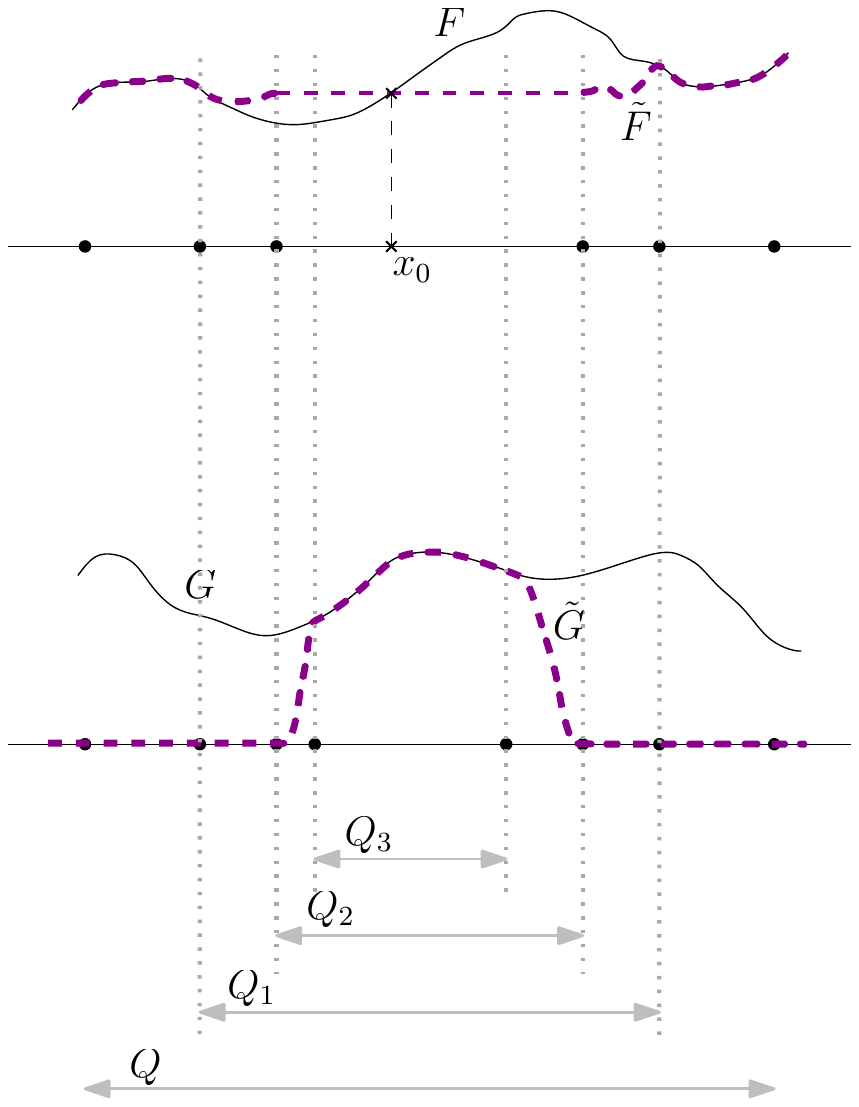}
	\caption{Producing $\til F$ and $\til G$ (the dashed lines).}
	\label{fig: construction_F_G_tildes_PB_flexibility}
\end{figure}

\subsubsection*{Construction of $\til{G}$}

Consider again a simplex $Q$ from our triangulation with subsets as mentioned,
$Q_3 \Subset Q_2 \Subset Q_1 \Subset Q$, satisfying $\Vol (Q \setm Q_3) \leq \Vol(Q) \cdot \eps$. Take a smooth function 
$\til{G} : Q \to \R$, $\til{G} = \psi \cdot G$, where
$\psi:Q\to [0,1]$
is a smooth function satisfying 
$\psi \restr _{Q_3} \equiv 1$, $\psi \restr _{Q \setm Q_2} = 0$.
Thus, we have 
$\til{G} \restr_{Q_3} \equiv G \restr_{Q_3}$, $\til{G} \restr _{Q \setm Q_2} \equiv 0$ and $|\til{G}(x)|\leq |G(x)|\ \forall x\in Q$.

Glue together all these patches of $\til{G}$ to get a smooth compactly supported function on $M$. The gluing is possible, since near the boundaries of each simplex, all the $\til{G}$-s vanish.\\

On a single simplex $Q$ we have 
\begin{align*} 
\int_Q | \til{G} - G |^q \om^n 
& = 
\int_{Q\setm Q_3} | \til{G} - G |^q \om^n 
\leq \int_{Q \setm Q_3} |G|^q \om^n 								
\leq \|G\|_{\infty}^q \int_{Q\setm Q_3} \om^n = \\
&	= \|G\|_{\infty}^q \cdot \Vol(Q \setm Q_3)
\leq \|G\|_{\infty}^q \cdot \Vol(Q) \cdot \eps \;.
\end{align*}

\noindent
Hence on the whole $M$ we get the bound
$$ \| \til{G} - G\|_q ^q = \int_M | \til{G} - G |^q \om^n 
\leq 
\|G\|_{\infty}^q \cdot \Vol( \{Q : Q \cap \supp(G) \neq \emptyset\} ) \cdot \eps \;,
$$
which depends on the volume of the union of all simplices intersecting $\supp (G)$.

Thus, by taking the diameter of the triangulation, $\delta$, small enough, we are able to produce pairs of Poisson commuting functions $\til{F},\til{G} \in C^\infty _c (M)$, so that $\til{F}$ is close to $F$ in the $C^0$-topology, and $\til{G}$ is close to $G$ in the $L_q$-norm.
We do have $\pois{ \til{F},\til{G} } = 0$, as for each $Q$, when $\til{G}$ is non-zero, $\til{F}$ is constant. Indeed, on each simplex $Q$ with the subset $Q_2$ as constructed, $\supp \til{G} \subseteq Q_2$ and $\til{F} \restr _{Q_2}$ is constant.

\end{proof} 









\section{$pb_4^q$: the two dimensional case} \label{sec: pb_4_q_two_dim}



In this section we prove \cref{thm: intro_thm_pb_4_dim=2_all} of the introduction.\\
Let $(M,\om)$ be a symplectic surface. We shall examine $pb_4^q$ of subsets inside $M$ both when $M$ has finite and infinite area.

Recall that for us, a curvilinear quadrilateral in $M$ is an image of a square $\hat{\Pi}$ by an embedding $\vphi: \hat{\Pi} \into M$.
Suppose that $X_0,X_1,Y_0,Y_1$ are sides of a curvilinear quadrilateral $\Pi \subset M$ (listed in counterclockwise order, here $\de \Pi = X_0 \cup Y_0 \cup X_1 \cup Y_1$).
We would like to show that $pb_4^q(\Pi):=pb_4^q(X_0,X_1,Y_0,Y_1)$ does not vanish and to compute it. We consider the cases $q=1$ and $1 < q < \infty$ separately at first, as we would use the result about upper bound for $q=1$ while proving the upper bound for $1<q<\infty$.\\

\begin{rmk} \label{rmk: def_F_4^q'}
	Recall that in the definition of $pb_4^q (\Pi)$ the infimum of $\| \pois{F,G} \|_q$ was taken over the set $\calF_4(\Pi) = \calF_4(X_0,X_1,Y_0,Y_1)$ of all pairs of functions 
	$(F,G)$ that satisfy 
	$$
	F \restr_{X_0} \leq 0,\ F \restr_{X_1} \geq 1,\ 
	G \restr_{Y_0} \leq 0,\ G \restr_{Y_1} \geq 1 \;.
	$$
	Instead, we can consider the infimum over a more restricted set, $\calF_4'(\Pi) = \calF_4'(X_0,X_1,Y_0,Y_1)$. This set consists of pairs $(F,G)$ of functions in $C^\infty_c (M)$, such that
	$0 \leq F,G\leq 1$ and 
	\begin{equation}
	F \restr _{near\ X_0} = G \restr _{near\ Y_0} = 0, 
	F \restr _{near\ X_1} = G \restr _{near\ Y_1} = 1 \;,
	\end{equation}
	
	where by saying "near" we mean in some neighborhood of the set.
	We will sometimes write $\calF_4' (\Pi,M)$ to emphasize that it is the set $\calF_4'(\Pi)$ with respect to $M$, i.e. that the functions $F$ and $G$ have compact support in $M$.
	
	\noindent
	We get an equivalent definition of $pb_4^q$ that is sometimes more convenient to use.
	The equivalence between these definitions can be proven repeating verbatim the proof in \cite[section 7.1]{polterov_rosen2014func_theory} (where it is given for the $L_\infty$-norm).\\
\end{rmk}

The following lemmas would be of use for us in the two-dimensional case.

%
%
%
%

\begin{lem} \label{lem: lem_lower_bd_L_1_pb_surface_wo_bdry}
	Let $(M,\om)$ be a symplectic surface without boundary, of finite or infinite area, and $\Pi \subset M$ a curvilinear quadrilateral with sides $X_0,X_1,Y_0,Y_1$ in counter-clockwise order. Then for any
	$(F,G) \in \calF_4'(\Pi) := \calF_4'(X_0,X_1,Y_0,Y_1)$ and $U$ being either $\Pi$ or $M\setminus \Pi$ we have 
	$$
	\int_U |\pois{F,G}|\om \geq 1 \;.
	$$
\end{lem}

\begin{proof} 
	
	Let $(F,G) \in \calF '_4(\Pi) := \calF_4' (X_0,X_1,Y_0,Y_1)$ be a pair of functions compactly supported in $M$.
	%
	By \cref{poisson_brck_def_1} for $n=1$, we have
	$ dF \wedge dG = -\pois{F,G} \om $. \\
	Using Stokes theorem and taking into account that for both options of $U$, 
	$\de U = \de(\Pi) = \de(M \setminus \Pi)$,
	$$
	\int_U |\pois{F,G}|\om \geq 
	\left| \int_U \pois{F,G} \om \right| = \left| \int_U dF \wedge dG \right|
	= \left| \int_U d(F dG) \right|
	= \left| \int_{\de U} FdG \right| 
	= \left| \int_{X_1} dG \right|
	= 1 \;.
	$$
\end{proof}

\begin{lem} \label{lem: lem_symp_stdrd_surface_w_curv_quadr} 
	
	Let $(M,\om)$ be a connected symplectic surface of area $B < \infty$, and let $\Pi \subset M$ be a closed curvilinear quadrilateral of area $A$. 
	Take any $A < C < B$ and an open rectangle $\Pi_C \subset \R^2$ of area $C$, with $\Pi_A \subset \Pi_C$ a closed rectangle of area $A$ (taking the standard area form in the plane).
	Then there exists an area preserving embedding 
	$\vphi : \Pi_C \to M$ such that $\vphi (\Pi_A) = \Pi$.
\end{lem}

The proof follows from Dacorogna-Moser theorem (see \cite{dacorogna_moser_1990}).

\begin{lem} \label{lem: symp_preserves_q_norm_of_PB_any_q_g_1}
	Let $(M, \om)$ be a symplectic surface and let $\Pi \subset M$ be a curvilinear quadrilateral.
	Take also $\Pi_0 \subset \R^2$ to be a closed square in the plane.
	Suppose that there exists a symplectomorphism $\vphi: U_{\Pi_0} \to U_{\Pi}$ from a neighborhood of $\Pi_0$ to a neighborhood of $\Pi$ in $M$ such that $\vphi (\Pi_0) = \Pi$.
	Let $(F,G) \in \calF_4(\Pi_0, U_{\Pi_0})$ (i.e., supported in $U_{\Pi_0}$), and define 
	$\til F = F \circ \vphi^{-1}$,
	$\til G = G \circ \vphi^{-1}$.
	Then $(\til F, \til G) \in \calF_4 (\Pi, U_{\Pi})$ (supported in $U_{\Pi}$) and
	$\| \pois{\til F,\til G} \|_q = \| \pois{F,G} \|_q$
	for any $1 \leq q < \infty$.
\end{lem}

\begin{proof}
	Denote $\psi = \vphi^{-1}$ and let $\om_{std}$ be the standard symplectic form on the plane.
	First, note that since $\vphi$ is a symplectomorphism, $\psi$ is also such, 
	therefore
	$\pois{\til F, \til G} (x) = \pois{F,G} (\psi(x))$.
	%
	%
	Hence, we have
	\begin{align*}
	\| \pois{\til F, \til G}\|_q^q 
	& =
	\int_M |\pois{\til F, \til G }|^q \om = 
	\int_{U_{\Pi}} |\pois{\til F, \til G}|^q \om
	=
	\int_{U_{\Pi}} | \psi^*  \pois{F,G}|^q \psi^*(\om_{std}) = \\
	& =
	\int_{U_\Pi} (\psi^*) (|\pois{F,G}|^q \om_{std})
	=
	\int_{U_{\Pi_0}} |\pois{F,G}|^q \om_{std} 
	=
	\| \pois{F,G} \|_q^q \;.		
	\end{align*}
	
\end{proof}


\subsection{$q = 1$}


\begin{thm} \label{thm: thm_pb_4_q=1_dim=2}
	For a symplectic surface $(M,\om)$ and a curvilinear quadrilateral $\Pi \subset M$, 
	$pb_4 ^{q=1}(\Pi) = 2 \;.$
\end{thm}

\begin{proof}[Proof of lower bound] 
	
	First, let us show that $2$ is a lower bound for $pb_4^{q=1} (\Pi)$. \\
	Take any $(F,G) \in \calF '_4(\Pi) := \calF_4' (X_0,X_1,Y_0,Y_1)$. Then by \cref{lem: lem_lower_bd_L_1_pb_surface_wo_bdry}
	$$
	\|\pois{F,G}\|_1 = \int_M |\pois{F,G}| \om = 
	\int_\Pi |\pois{F,G}| \om + \int_{M \setm \Pi} |\pois{F,G}| \om
	\geq 2 \;,
	$$
	
	\noindent
	therefore $pb_4 ^{1} (\Pi) \geq 2$.
\end{proof}

We would like to show that $2$ is also an upper bound for $pb_4 ^{1} (\Pi)$. 
The proof would be very similar to the proof of the upper bound in the case of $1 < q <\infty$ below. 
Therefore, we will first show the upper bound for $1 < q < \infty$ (see \cref{thm: thm_q>1_dim=2_pb_4_finite_area}) 
and then deduce the limiting case $q=1$ from the same construction.

\subsection{$1 < q < \infty$}

%

We study $pb_4^q (\Pi) $ of a curvilinear quadrilateral $\Pi \subseteq M$ on a connected surface without boundary for $1 < q < \infty$. In this case, $pb_4^q$ appears to depend on the areas of $\Pi$ and $M$.
We first consider the case when $\Area M < \infty$ and then use it for the case of a surface of infinite area (see \cref{thm: thm_pb_4_q>1_dim=2_infinite_area} below). 

\begin{thm} \label{thm: thm_q>1_dim=2_pb_4_finite_area}
	Let $1 < q < \infty$. Denote $A = \Area(\Pi)$ and $B = \Area(M) < \infty$.
	Then 
	\begin{align} \label{lower_bound_q_g_1_statement}
	\displaystyle{pb_4^q (\Pi) = \Big( \frac{1}{A^{q-1}} + \frac{1}{(B-A)^{q-1}} \Big)^{1/q}} \;.
	\end{align} 
\end{thm}

\begin{proof}
	
	First, we would show that the right-hand-side of \cref{lower_bound_q_g_1_statement} is a lower bound for $pb_4^q(\Pi)$. 
	Take any $(F,G)\in \calF_4'(\Pi)$, a pair of functions compactly supported inside $M$.
	By \cref{lem: lem_lower_bd_L_1_pb_surface_wo_bdry} applied to $U$ being either $\Pi$ or $M\setminus \Pi$, we have
	$
	\int_U |\pois{F,G}| \geq 1
	$.
	
	Let $p$ be such that $\frac{1}{q}+\frac{1}{p} = 1$ (then $q/p = 1-q$). 
	Let us note that for any smooth function $f$ on $U$, by H\"{o}lder inequality we have
	$$
	\int_U |f| \om \leq \Big( \int_U |f|^q  \om \Big)^{1/q} \cdot \Big( \int_U |1|^p  \om \Big)^{1/p}
	= \Big( \int_U |f|^q \om \Big)^{1/q} \cdot \Big( \Area(U) \Big)^{1/p} \;,
	$$
	so 
	$$\displaystyle{ 
		\Big( \int_U |f|^q \om \Big)^{1/q} \geq 
		\frac{\int_U |f| \om}{\Big( \Area(U) \Big)^{1/p}} 
	} \;.
	$$
	In our case, for $f=\pois{F,G}$ we get
	$$
	\int_U |\pois{F,G} |^q \om 
	\geq 
	\frac{ \Big( \int_U |\pois{F,G}| \om \Big)^q }{\Big( \Area(U) \Big)^{q/p}}
	\geq
	\frac{ 1 }{ \Big( \Area(U) \Big)^{q-1} } \;.
	$$
	
	\noindent
	Hence, 
	$$
	\int_{\Pi} |\pois{F,G} |^q \om \geq \frac{1}{ A^{q-1} }
	,\ \ 
	\int_{M\setm \Pi} |\pois{F,G} |^q \om \geq \frac{1}{ (B-A)^{q-1} } \;,
	$$
	and overall we have
	$$
	\| \pois{F,G}\|_q = 
	\Big( \int_M |\pois{F,G}|^q \om \Big)^{1/q} 
	\geq
	\Big( \frac{1}{ A^{q-1}} + \frac{1}{ (B-A)^{q-1}} \Big)^{1/q} \;.
	$$\\

	In order to prove that an equality in \cref{lower_bound_q_g_1_statement} holds, we shall construct pairs of functions $F,G \in C^\infty_c (M)$ with $\| \pois{F,G } \|_q$ arbitrary close to the declared value of $pb_4^q (\Pi)$.
	
	
	We first present a construction for a rectangle $\Pi$ of area $A$ contained in another rectangle $M$ in the plane of area $B$.\\

	Fix $\epsi>0$ and $A<C<B$. Let $\Pi$ be a rectangle in the plane, 
	$\Pi = [0,A]\times[0,1]$, and let $K$ be
	$\Pi \subset K=[-\epsi,C+\epsi]\times[-2\epsi, 1+2\epsi] \subset int(M)$.\\

	Define the following four smooth functions:
	\footnote{
		Following a construction by Lev Buhovsky as presented in \cite[section 7.5.3]{polterov_rosen2014func_theory}.
	} 
	\begin{itemize}
		\item 
		$u_1 : \R \to [0,1]$, such that 
		$\supp (u_1) \subset (0,C)$, 
		$u_1 (A) = 1$. Later, a more specific function with this properties will be considered.
		\item
		$v_1 : \R \to [0,1]$, such that
		$\supp (v_1) \subset (-\eps, 1+\eps)$ and
		$v_1 \restr_{[0,1]} \equiv 1$.
		\item
		$u_2 : \R \to [0,1]$, such that
		$\supp (u_2) \subset (-\eps, C + \epsi)$ and
		$u_2 \restr_{[0,C]} \equiv 1$.
		\item
		$v_2 : \R \to [-\eps, 1+\eps]$, such that
		$\supp (v_2) \subset (-2\eps, 1+2\eps)$ and
		$v_2 \restr_{[-\eps,1+\eps]} = id$.
	\end{itemize}
	
	\begin{figure} [h!1]
		\centering
		\includegraphics[scale = 1, trim = 0 0 0 0 ]{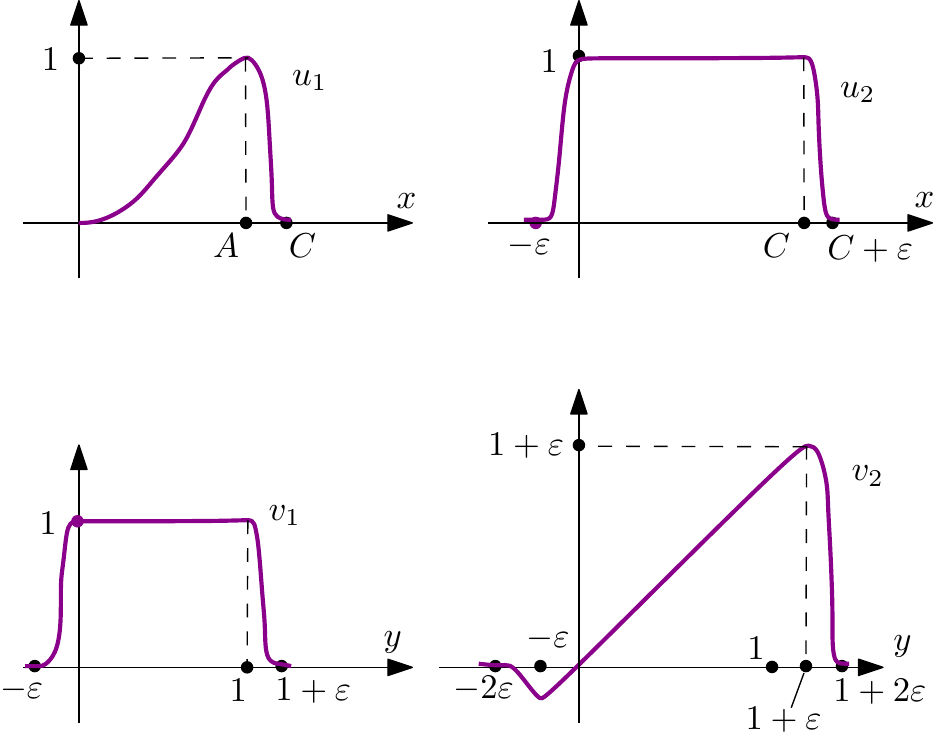}
		\caption{Constructing $F,G$ using four functions $u_1, v_1, u_2, v_2$.}
		\label{fig: four_functions_rectangles.}
	\end{figure}
	
	Put $F(x,y) = u_1(x) v_1(y)$, $G(x,y) = u_2(x) v_2(y)$. These functions belong to $C^\infty_c(M)$, they are supported in $K$, and $(F,G) \in \calF_4'(\Pi)$ (note that $K$ depends on $\epsi$ and $C$).
	
	\noindent
	We have
	$$
	\pois{F,G} = \underbrace{u_1(x) u_2'(x)}_{= 0} v_1'(y) v_2(y) 
	- \underbrace{u_1'(x) u_2(x)}_{=u_1'(x)} \underbrace{v_1 (y) v_2'(y)}_{=v_1(y)} 
	= - u_1'(x) v_1(y) \;.
	$$
	Hence 
	\begin{align*}
	\| \pois{F,G} \|_q^q &= \int_M |\pois{F,G}|^q = \int_K |u_1'(x)|^q \cdot |v_1(y)|^q dxdy = \\
	&= \int_0^C |u_1'|^q \cdot \int_{-\epsi}^{1+\epsi} |v_1|^q 
	\leq (1+2\epsi) \int_0^C |u_1'|^q	 \;.
	\end{align*}
	
	Observe that if we take $u_1$ to be linear on $[0,A]$ and on $[A,C]$, i.e. increasing from 0 to 1 on $[0,A]$ and decreasing back to zero at $C$, we would get
	$\int_0^C |u_1'|^q = \int_0^A |u_1'|^q  +  \int_A^C |u_1'|^q = \frac{1}{A^{q-1}} + \frac{1}{(C-A)^{q-1}}$, 
	and then 
	$\| \pois{F,G} \|_q^q \leq \Big( \frac{1}{A^{q-1}} + \frac{1}{(C-A)^{q-1}} \Big)(1+2\epsi)$. 
	We can approximate this piece-wise linear $u_1$ in the $L_\infty$-topology by smooth functions to obtain 
	$\int_0^C |u_1'|^q $ arbitrarily close to $\frac{1}{A^{q-1}} + \frac{1}{(C-A)^{q-1}}$.\\
	For instance, take $u_1$ to be 
	linear on $[2\eps, A- 2\eps]$ with 
	$u_1 (2\eps) = \eps$ and $u_1(A-2\eps) = 1-\eps$,
	and 
	linear on $[A+2\eps, C-2\eps]$ with 
	$u_1(A+2\eps)=1-\eps$ and $u_1(C-2\eps)=\eps$.
	Then smoothly extend it to an increasing function on the whole $[0,C]$, such that $u_1' \leq 1$ and $u_1' = 0$ close to $0$, $A$ and $C$ (taking $\eps$ small enough with respect to $A$ and $C$).\\
	The slopes on the linear parts would then be 
	$m_{[2\eps,A-2\eps]} = \frac{1 - 2\eps}{A - 4\eps}$ and
	$m_{[A+2\eps,C-2\eps]} = -\frac{1 - 2\eps}{(C-A) - 4\eps}$.
	And hence
	\begin{align*}
	\int_0^C |u_1'|^q & \leq 
	8\epsi \cdot 1 
	+ (A-4\eps)\cdot \frac{(1-2\eps)^q}{(A-4\eps)^q}
	+ ((C-A)-4\eps)\cdot \frac{(1-2\eps)^q}{((C-A)-4\eps)^q} \\
	& = 
	8\epsi
	+ (1-2\eps)^q \cdot 
	\left(
	\frac{1}{(A-4\eps)^{q-1}}
	+ \frac{1}{((C-A)-4\eps)^{q-1}}
	\right) \;.
	\end{align*}

	Thus, for any $1 < q < \infty$, taking $\eps \to 0$ and $C\to B$, 
	we would get pairs $(F,G) \in \calF_4'(\Pi)$ with 
	$ \| \pois{F,G}\|_q $ arbitrarily close to
	$
	\Big( \frac{1}{A^{q-1}} + \frac{1}{(B-A)^{q-1}} \Big)^{1/q} \;.\\ \\
	$

	This proves \cref{thm: thm_q>1_dim=2_pb_4_finite_area} for $1 < q < \infty$ and for this model of rectangle inside another rectangle in the plane.\\

	Let us go back to the general case. 
	We have a symplectic surface $(M,\om)$ without boundary of finite area $B$ and a curvilinear quadrilateral $\Pi\subseteq M$ of area $A$.

	\noindent
	Take any $A < C <B$ and consider an open rectangle $\Pi_C \subset \R^2$ of Euclidean area $C$.
	By \cref{lem: lem_symp_stdrd_surface_w_curv_quadr}, 
	there exists an area preserving map $\vphi : \Pi_C \to M$ that takes a rectangle $\Pi_A$ of area $A$ to $\Pi$. 
	Note that the map $\vphi: \Pi_C \to \vphi (\Pi_C) $ is a symplectomorphism.
	%
	%
	
	If $ (F,G) \in \calF_4'(\Pi_A, \Pi_C)$ (i.e. supported inside $\Pi_C$, see \cref{rmk: def_F_4^q'}) 
	, take
	$\til F = F \circ \vphi$ and
	$\til G = G \circ \vphi$. 
	Then $ (\til F, \til G) \in \calF_4' (\Pi, \vphi(\Pi_C))$.
	
	\noindent
	Using \cref{lem: symp_preserves_q_norm_of_PB_any_q_g_1} we can conclude that
	$$
	pb_4^q(\Pi,M) \leq 
	pb_4^q (\Pi, \vphi (\Pi_C)) = pb_4^q (\Pi_A, \Pi_C) = \left( \frac{1}{A^{q-1}} + \frac{1}{(C-A)^{q-1}} \right)^{1/q} \;.
	$$
	
	\noindent
	Therefore, taking $C \to B$ we get that
	$pb_4^q (\Pi, M) \leq \left( \frac{1}{A^{q-1}} + \frac{1}{(B-A)^{q-1}} \right)^{1/q}$. But we have already shown the opposite inequality, hence we have the equality \cref{lower_bound_q_g_1_statement}.\\

	To get the upper bound $2$ for $q=1$ we can apply the same construction (putting $q=1$ everywhere), both for the special case of rectangles in the plane and for the general case.
	
\end{proof} 

\begin{rmk}
	For $q=1$, it is enough to have a diffeomorphism $\vphi: \Pi_B \to M$ with the above properties, instead of a symplectomorphism, as the statement of \cref{lem: symp_preserves_q_norm_of_PB_any_q_g_1} would hold for a diffeomorphism in this case.
\end{rmk}

As a corollary of \cref{thm: thm_q>1_dim=2_pb_4_finite_area}, we will be able to compute $pb_4^q$ for the case of a surface with infinite area.

\begin{thm} \label{thm: thm_pb_4_q>1_dim=2_infinite_area}
	Suppose $(M, \om)$ is a connected symplectic surface without boundary of infinite area.
	Then for a curvilinear quadrilateral $\Pi \subset M$ of area $A$,
	for $1 < q < \infty$, 
	${pb_4^q(\Pi) =  \left( \frac{1}{A^{q-1}} \right)^{1/q} }\ $.
\end{thm}

\begin{proof}
	Let 
	\footnote{
		Such $M_1$ exists for instance by \cref{lem: lem_symp_stdrd_surface_w_curv_quadr}.
	}
	$M_1 \subset M$ be a connected subsurface with finite area $B > A$, so that $\Pi \subset int(M_1)$.
	Denote by $pb_4^q(\Pi,M)$ this invariant with respect to functions that have compact support inside $M$, and similarly $pb_4^q(\Pi, M_1)$ for $M_1$. 
	
	Since $M_1 \subset M$ and by \cref{thm: thm_q>1_dim=2_pb_4_finite_area}, we have
	$$
	pb_4^q (\Pi,M) \leq p_4^q (\Pi, M_1) = 
	\left( \frac{1}{A^{q-1}} + \frac{1}{(B-A)^{q-1}} \right)^{1/q} \;.
	$$
	This holds for any $B > A$, hence we get an upper bound on $pb_4^q (\Pi,M)$,
	$$
	pb_4^q (\Pi,M) \leq  
	\left( \frac{1}{A^{q-1}} \right)^{1/q} \;.
	$$
	
	Let us now show that actually an equality holds in this last inequality. 
	Suppose on the contrary that 
	$pb_4^q (\Pi,M) = \left( \frac{1}{A^{q-1}} \right)^{1/q} - \epsi$
	for some $\epsi > 0$.
	Then there exist two functions $(F,G) \in \calF_4 (M)$ with
	$
	\|\pois{F,G}\|_q \leq \left( \frac{1}{A^{q-1}} \right)^{1/q} - \frac{\epsi}{2} \;.
	$ 
	
	\noindent	
	Observe that $F$ and $G$ have compact support in $M$. 
	Consider some $M_1$ of area $B$ diffeomorphic to an open disk, such that
	$ M_1 \supset \supp F \cup  \supp G \cup \Pi \supset \Pi$.
	Since $(F,G) \in \calF_4 (\Pi, M_1)$, by \cref{thm: thm_q>1_dim=2_pb_4_finite_area} we get
	$$
	pb_4^q (\Pi, M_1) \leq \| \pois{F,G} \|_q 
	\leq  \left( \frac{1}{A^{q-1}} \right)^{1/q} - \frac{\epsi}{2} 
	<  \left( \frac{1}{A^{q-1}} + \frac{1}{(B-A)^{q-1}} \right)^{1/q}
	= pb_4^q (\Pi, M_1) \;,
	$$
	which is a contradiction.
	Hence $pb_4^q (\Pi,M) = \left( \frac{1}{A^{q-1}} \right)^{1/q}$.
\end{proof}

\begin{rmk}
	For $q=1$ and $(M,\om)$ of infinite area, we can use the same proof to obtain that 
	$pb_4^1 (\Pi,M) = 2$.
\end{rmk}

\begin{rmk}
	Consider a fixed curvilinear quadrilateral $\Pi$ on a symplectic surface $M$, still in the setting of \cref{thm: thm_q>1_dim=2_pb_4_finite_area}.
	Let us note that using the values computed for $pb_4^q(\Pi)$ we can find a lower bound on $pb_4^\infty(\Pi)$. More precisely, the following inequality holds:
	\begin{equation}
	\label{eq: lower_bound_of_pb_4^infty}
	pb_4^\infty (\Pi) \geq \limsup_{q\to \infty} pb_4^q (\Pi)
	\end{equation}
	
	\noindent
	Indeed, take any $(F,G) \in \calF_4'(\Pi)$. Then by definition
	$
	\| \pois{F,G} \|_q \geq pb_4^q(\Pi)
	$.
	Taking $q\to \infty$ we have, for fixed $(F,G)$,
	$
	\| \pois{F,G} \|_\infty \geq \limsup_{q\to \infty} (pb_4^q (\Pi))
	$.
	This is true for any $(F,G)$, 
	hence \cref{eq: lower_bound_of_pb_4^infty} holds.
	
	In our case, this gives the following precise lower bound on $pb_4^\infty(\Pi)$, which was already proven (see e.g. \cite[7.5.3]{polterov_rosen2014func_theory}):
	\begin{equation*}
	\label{eq:lower_bound_pb_4^infty_conclusion}
	pb_4^\infty (\Pi) \geq  \limsup_{q\to \infty} (pb_4^q (\Pi))
	= \lim_{q\to \infty} \Big( \frac{1}{A^{q-1}} + \frac{1}{(B-A)^{q-1}} \Big)^{1/q} 
	= \max \left( \frac{1}{A}, \frac{1}{B-A} \right)
	\;.
	\end{equation*}
	
	Similarly, we observe that the function $q\mapsto pb_4^q(\Pi)$ is upper semi-continuous.
\end{rmk}







\section {$pb_4^q$: the multi-dimensional case}

Consider a symplectic manifold $(M^{2n},\om)$, where $n\geq 2$. 
Let $X_0, X_1, Y_0, Y_1$ be compact subsets of $M$, such that 
$X_0 \cap X_1 = Y_0 \cap Y_1 = \emptyset$, assuming also that $X_1$ is a submanifold with or without boundary.
%

\begin{thm} \label{thm: thm_pb_4_vanish_high_dim}
	Denote $d = \dim X_1$ and assume that $d \leq 2n-2$. Then $pb_4^q(X_0, X_1, Y_0, Y_1) = 0$ whenever $q \leq 2n - d, n\geq 2$.
\end{thm}

Let $G\in C_c^{\infty}(M)$ be any function that assumes values in $[0,1]$ such that
$G\restr_{near\ Y_0} \equiv 0$ and $G\restr_{near\ Y_1} \equiv 1$.
To prove \cref{thm: thm_pb_4_vanish_high_dim}, 
we show below that there exists $F\in C^\infty_c(M)$ so that $(F,G)\in \calF_4'(X_0,X_1,Y_0,Y_1)$ with arbitrarily small $\|\pois{F,G}\|_q$, concluding that $pb_4^q (X_0,X_1,Y_0,Y_1)$ vanishes. \\
%
%

We start with a few lemmas.

\begin{lem}\label{lem_func_fast_decay_dim1}
	Fix $2\leq m \in \N$. 
	Then for any $\eps, \delta >0$ and $1\leq k \leq m$ there exists a non-negative function $f\in C^{\infty} ([0,+\infty))$ supported in $[0,\delta)$ with $\max |f| = f(0) = 1$, such that 
	$\int _0 ^{\infty} |f'(r)|^k r^{m-1} \leq \eps$,
	$\int _0 ^{\infty} |f|^k \leq \epsi$,
	and such that $f$ is constant near $0$.
\end{lem}

\begin{proof} 
	Take any smooth function $ h : [0,\infty) \rightarrow [0,\infty) $, such that $ h = 1 $ near $ 0 $ and such that the support of $ h $ is contained in $ [0,1/2] $. 
	For every $ \alpha > 0 $ define $ h_\alpha(r) = h(r^\alpha) $. 
	Since the support of $ h_\alpha $ lies in $ \left[0,\frac{1}{2^{1/\alpha}}\right] $
	we have
	$$ 
	\int_0^\infty h_\alpha(r)^k \, dr \leqslant \frac{\|h_\alpha^k \|_\infty}{2^{1/\alpha}} = \frac{\|h^k\|_\infty}{2^{1/\alpha}} \;.
	$$ 
	In particular, the left-hand-side of the inequality tends to zero as 
	$ \alpha \rightarrow 0+ $. 
	We also have 
	$$ \int_0^\infty | h_\alpha'(r)|^k r^{m-1} \, dr = 
	\alpha^k \int_0^\infty | h'(r^\alpha)|^k r^{m + k\alpha- k - 1} \, dr =  
	$$
	$$ = \alpha^{k-1} \int_0^\infty | h'(t) |^k t^{\frac{m-k}{\alpha} +k-1} \, dt  \leqslant \frac{\alpha^{k-1} \|h'\|_\infty^k}{2^{\frac{m-k}{\alpha} + k - 1}} \;,
	$$ 
	where in the second equality we made the substitution $ t = r^\alpha $, 
	and in the last step we estimated the integral from above by the maximum of the integrand,
	taking into account that $t \in [0,\frac{1}{2}]$. 
	For $ m > 1 $ and $ 1 \leqslant k \leqslant m $, the upper bound converges to $ 0 $ 
	as $ \alpha \rightarrow 0+ $. 
	Hence, the integral 
	$ \int_0^\infty | h_\alpha'(r)|^k r^{m-1} \, dr $ 
	converges to $ 0 $ as well, when $ \alpha \rightarrow 0+ $. 
	We conclude that $ f := h_\alpha $ will satisfy all the requirements, 
	taking small enough $ \alpha > 0 $.
\end{proof}

The next lemma is, in a sense, a generalization of the previous one to higher dimensions.

\begin{lem} \label{lem_func_fast_decay_high_dim}
	Fix $2\leq m \in \N$. For all $\eps, \delta >0$ and $1\leq k \leq m$ there exists a non-negative function $f\in C^\infty (\R ^m)$ supported in a ball $B_\delta$ of radius $\delta$ around 0, with $\max |f| = f(0) = 1$, such that 
	$\int_{\R^m} \|\na f\|^k d\Vol \leq \eps$
	and
	$\int_{\R ^m} |f|^k dVol \leq \epsi$
	.
\end{lem}

\begin{proof}
	%
	Let $\epsi, \delta > 0$.
	Consider $\R ^m$ with polar coordinates, consisting of the radius $R(x) = \|x\|$ and coordinates on the unit sphere $S^{m-1}$, so that the volume element in $\R ^m$ rewrites as $d\Vol = r^{m-1} dr d\sigma$, where $d\sigma$ is the volume element on $S^{m-1}$.\\
	
	Let $g=g(r)$ be a function that satisfies the requirements of 
	\cref{lem_func_fast_decay_dim1} for our $\epsi, \delta >0$.
	Take the radial function $f: \R ^m \to \R$ defined by $f(x) = g(r(x))$.\\
	Then $f$ is a smooth function on $\R ^m$, supported in $B_\delta$, with $\max |f| = f(0) = 1$. Let us verify that it also satisfies the other two declared properties.\\
	%
	%
	Note that for a radial function we have $\na f = \frac{\de f}{\de r} \vect{e_r}$, where $\vect{e_r}$ is the unit vector in the radial direction. And so, $\|\na f\| = |\frac{\de f}{\de r}| = |g'(r)|$.
	%
	We consider the integral 
	\begin{align*}
	\int_{\R ^m} \|\na f\|^k d\Vol 
	& = 
	\int_0 ^\infty  \int_{S^{m-1}} \|\na f\|^k r^{m-1} d\sigma dr 
	=
	\int_{S^{m-1}} d\sigma \int_0^\infty \|\na f\|^k r^{m-1} dr = \\
	&
	= C_m \int_0^\infty | g'(r) |^k r^{m-1} dr \;,
	\end{align*}
	
	\noindent
	where $C_{m} = vol(S^{m-1})$ is a constant independent of $k$ and $f$.
	The right-hand-side can be made as small as needed, by \cref{lem_func_fast_decay_dim1} $\forall 1\leq k\leq m$.
	\\
	
	Finally, for the integral $\int_{\R ^m} |f|^k dVol$, 
	calculating again in polar coordinates, we have
	\begin{align*}
	\int_{\R ^m} |f| ^k d\Vol &=
	\int_0^{\infty} \int_{S^{m-1}} |f(r)|^k r^{m-1} d\sigma dr 
	= \int_{S^{m-1}} d\sigma \int_0^\infty |g(r)|^k r^{m-1} dr = \\
	& = C_m \int_0^{\delta} |g(r)|^k r^{m-1} dr \leq C_m \cdot \delta^{m-1} \int_0^\delta |g(r)|^k dr
	\;,
	\end{align*}
	and the expression in the right-hand-side can be made arbitrary small, by \cref{lem_func_fast_decay_dim1}.
	
\end{proof}

\begin{rmk} \label{equivalence_of_metrics}
	Since on any compact set $B \subset \R^m$ any two Riemannian metrics are equivalent, the statement of \cref{lem_func_fast_decay_high_dim} holds true not only for the Euclidean metric, but for any other Riemannian metric on $B$.
\end{rmk}


\begin{rmk}
	At the beginning of the proof of \cref{thm: thm_pb_4_vanish_high_dim}, we will use the following basic notion.
	Let $(M,\om)$ be a symplectic manifold.
	Having a function $F\in C^\infty(M)$, we can define a smooth vector field on $M$ associated with $F$. We consider a vector field $\sgrad F$ that satisfies the identity
	$$
	\om (\sgrad F, \cdot) = -dF(\cdot)\;.
	$$
	Such a vector field exists and it is unique, by the non-degeneracy of $\om$.
	It is called the Hamiltonian vector field of $F$.
\end{rmk}

Let us mention here that $m \in \N$ as appears in \cref{lem_func_fast_decay_high_dim} will play the role of 
$\codim X_1 = 2n-d$ in the following proof.
\label{proof: prf_thm_pb_4_vanish_high_dim}
\begin{proof} [Proof of \cref{thm: thm_pb_4_vanish_high_dim}]
	
	Our general strategy will be as follows. For any $G\in C^{\infty}_c (M)$ such that $0\leq G \leq 1, G\restr_{near\ Y_0} = 0, G \restr_{near \ Y_1}=1$, 
	we want to find a function $F\in C^{\infty}_c(M)$ so that 
	$(F,G) \in \calF_4'(\Pi)$
	and
	$\|\pois{F,G}\|_q$ is arbitrarily small.\\
	%
	%
	On $M$, pick a Riemannian metric $\rho$. We consider the norm $\| \cdot \|_\rho$ and the gradient $\nabla_\rho$ with respect to this metric.\\
	
	%
	%
	
	Note that by the definition of $\nabla_\rho$ and by the Cauchy-Schwartz inequality,
	\begin{align}
	\begin{split}
	\label{eq: bounding_PB_using_gradient}
	| \pois{F,G} | &= |\om(\sgrad F,\sgrad G)| = | dF(\sgrad G) | = \\ 
	& = | (\nabla_\rho F, \sgrad G)_\rho | 
	\leq \| \nabla_\rho F \|_\rho \cdot \| \sgrad G \|_\rho \;.
	\end{split}
	\end{align}
	
	\noindent
	Hence,
	\begin{align}
	\| \pois{F,G} \|_q^q = \int_M | \pois{F,G} |^q \om^n
	\leq \max_M \| \sgrad G \| \cdot \int_M \| \nabla_\rho F \|^q \om^n \;,
	\end{align}
	
	\noindent
	and it would be enough to produce a function $F$ as above with arbitrary small 
	$\| \nabla_\rho F\|_q$. 
	We will do so by first constructing appropriate functions locally in a neighborhood of $X_1$ (using \cref{lem_func_fast_decay_high_dim}), and then gluing them.\\

	Cover $X_1$ by a finite collection $\{U_\alpha\}$ of open subsets of $M$, each equipped with a diffeomorphism $\ph_\alpha:U_\alpha \to \R^{2n}$ that flatten $X_1$ in the following sense.
	Take coordinates $z_1,\ldots, z_d, z_{d+1}, \ldots, z_{2n}$ with respect to the standard basis $e_1, \ldots , e_{2n}$ on $\R^{2n}$. 
	We require $\vphi_\alpha$ to satisfy
	$\vphi_\alpha (X_1 \cap U_\alpha) \subset \{z_{d+1} = \ldots = z_{2n} = 0\}$. 
	Suppose also that the sets $U_\alpha$ are all small enough so that $U_\alpha \cap X_0 = \emptyset \ \forall \alpha$.
	(See \cref{fig: cover_flattening_projection}.)
	
	Take a collection of cutoff functions $\{\eta_\alpha : U_\alpha \to [0,1] \}$ that form a partition of unity subordinate to the cover $\{U_\alpha\}$ of $X_1$, so that  
	$\supp \eta_\alpha \subset U_\alpha$ and 
	$\forall x\in X_1, \  \sum_{\alpha} \eta_\alpha (x) = 1$. 
	
	Let us emphasize in advance that the cover $\{ (U_\alpha, \vphi_\alpha) \}$ and the collection $\{\eta_\alpha\}$ are fixed throughout the proof.\\
	
	We want to construct suitable functions on each $U_\alpha$ separately, and then glue them to a function $F:M \to \R$, using $\{ \eta_\alpha \}$.

	\begin{figure} [h!]
		\centering
		\includegraphics[scale= 0.75]{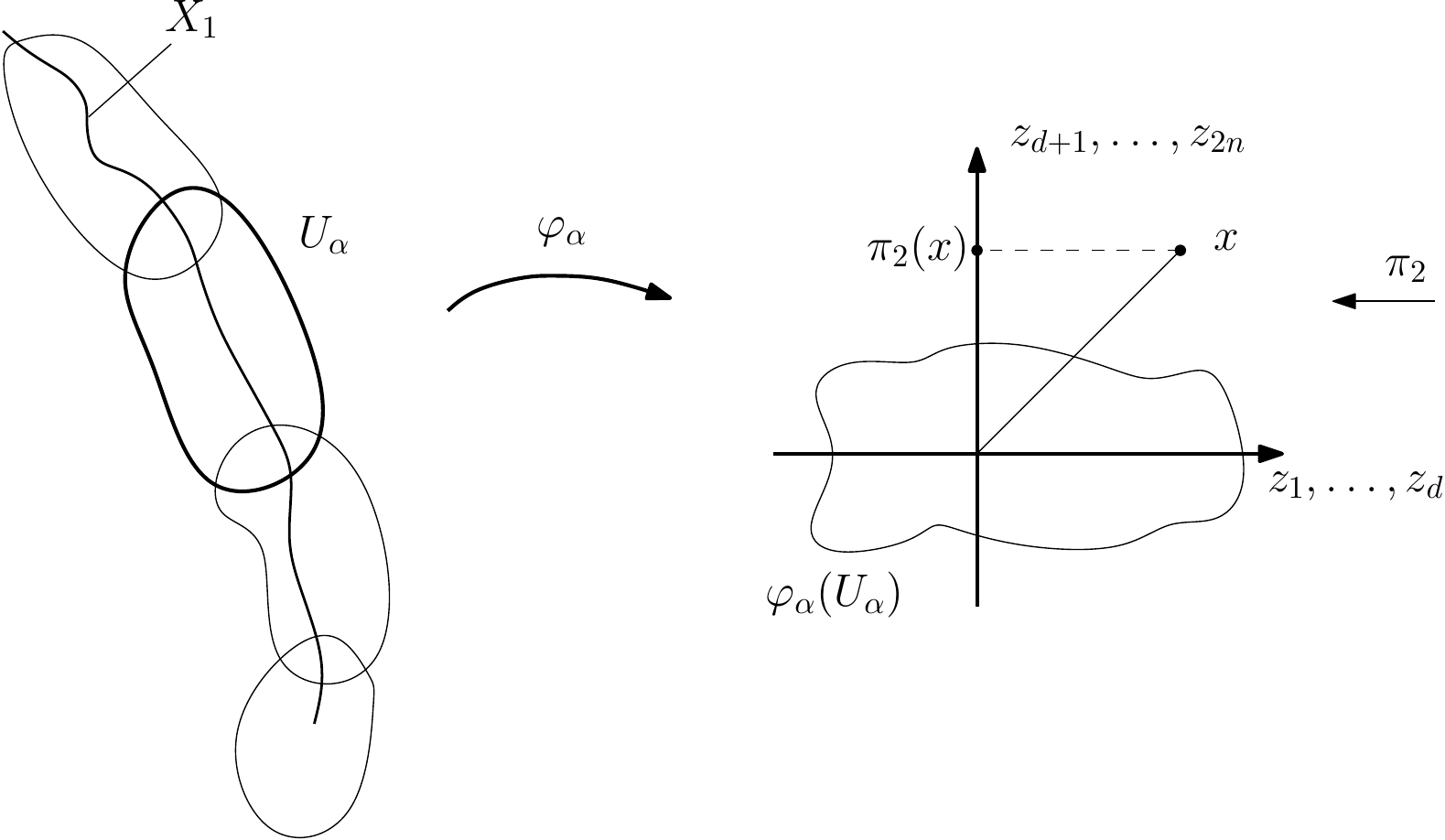}
		\caption{Covering $X_1$ by $\{ U_\alpha, \vphi_\alpha \}$ and constructing $F$ locally.}
		\label{fig: cover_flattening_projection}
	\end{figure}

	\noindent
	Consider a single $\vphi_\alpha(U_\alpha)$. 
	%
	Let $\pi_2: \vphi_\alpha (U_\alpha) \to span\{e_{d+1},\ldots, e_{2n}\}$ be the projection to $span\{e_{d+1},\ldots, e_{2n}\}$.
	Denote $r(x) = \| \pi_2(x) \| = \sqrt{z_{d+1}^2(x) + \ldots + z_{2n}^2(x)}$.\\ 
	
	The function we want to define on $\vphi_\alpha (U_\alpha)$ would depend only on the distance of a point $x$ from $span\{e_{d+1},\ldots, e_{2n}\}$, i.e., on $r(x)$.
	Let 
	$F _\alpha : \vphi_\alpha (U_\alpha) \to \R$ 
	be defined by
	$F _\alpha (x) = f_\alpha (\pi_2(x)) = f_\alpha ( z_{d+1}(x), \ldots, z_{2n}(x) )$,
	where $f_\alpha$ is a function that fulfills the requirements of \cref{lem_func_fast_decay_high_dim} \ 
	for $\epsi>0$ and $\delta = \max_{\vphi_\alpha(U_\alpha)} \{r(x)\} $, and the linear space $span\{e_{d+1},\ldots, e_{2n}\}$, i.e. for $m=2n-d$ in the notations of the lemma. 
	(In fact, the value of $\delta$ is not important for the construction.)\\
	%
	%
	%


	We look now at $U_\alpha$ and take the pullback of $F_\alpha$,
	$\til F_\alpha = F_\alpha \circ \vphi_\alpha$, thus defining 
	$\til F_\alpha : U_\alpha \to \R$. 
	Extend it by zero to the whole $M$. 

	Consider $F : M \to \R$ defined by $F = \sum_{\alpha} \eta_\alpha \til F_\alpha$. \\
	Observe the following properties of~$F$. 
	It is a smooth function with compact support that is contained in $\cup U_\alpha$. 
	Also, 
	$F \restr_{X_0} \equiv 0$, since 
	$X_0 \cap (\cup U_\alpha) = \emptyset$, 
	and
	$F \restr_{X_1} \equiv 1$
	since 
	$ F_\alpha \restr_{\vphi_\alpha (X_1 \cap U_\alpha)} \equiv 1$ and 
	$\sum \eta_\alpha \equiv 1$ on $X_1$.
	Thus, $(F,G) \in \calF_4' (\Pi) = \calF_4' (X_0,X_1,Y_0,Y_1)$.\\

	%
	

	We have 
	$
	\na_\rho F = \na_\rho ( \sum \eta_\alpha \til F_\alpha ) 
	= \sum \na_\rho \eta_\alpha \cdot \til F_\alpha + \sum \eta_\alpha \na_\rho \til F_\alpha
	$, and by the triangle inequality,
	$$
	{\| \na_\rho F \|_\rho \leq  
		\sum_\alpha \| \na_\rho \eta_\alpha \cdot \til F_\alpha \|_\rho + 
		\sum_\alpha \| \eta_\alpha \na_\rho \til F_\alpha\|_\rho } \;.
	$$
	
	%
	
	Combining this with the Minkowski inequality and then using the positivity of the integrands, we have
	
	\begin{align*}
	\left( \int_{M} \| \na_\rho F \|_\rho^q \om^n \right)^{1/q} & 
	\leq 
	\left( \int_{\cup U_\beta} 
	\left( \sum_\alpha \| \na_\rho \eta_\alpha \til F_\alpha \|_\rho + 
	\sum_\alpha \| \eta_\alpha \na_\rho \til F_\alpha \|_\rho  \right)^q \om^n 
	\right)^{1/q} \leq \\
	& \leq
	%
	%
	%
	%
	\sum_{\alpha} \left( \int_{\cup U_\beta} 
	\| \na_\rho \eta_\alpha \cdot \til F_\alpha \|_\rho^q \ \om^n \right)^{1/q}
	+
	\sum_{\alpha} \left( \int_{\cup U_\beta} 
	\| \eta_\alpha \na_\rho \til F_\alpha \|_\rho^q	\ \om^n \right)^{1/q} = \\
	&
	\leq
	\sum_{\alpha} \left( \int_{U_\beta} 
	\| \na_\rho \eta_\alpha \cdot \til F_\alpha \|_\rho^q \ \om^n\right)^{1/q}
	+
	\sum_{\alpha} \left( \int_{U_\beta} 
	\| \eta_\alpha \na_\rho \til F_\alpha \|_\rho^q	\ \om^n \right)^{1/q} = \\
	&
	= 
	\sum_{\alpha} \left( \int_{U_\beta}  |\til F_\alpha|^q 
	\| \na_\rho \eta_\alpha \|_\rho^q \ \om^n \right)^{1/q}
	+
	\sum_{\alpha} \left( \int_{U_\beta} |\eta_\alpha|^q 
	\| \na_\rho \til F_\alpha \|_\rho^q	\ \om^n \right)^{1/q} \;.
	\end{align*}
	
	Since there is a finite number of sets in the covering 
	and since $\{ \eta_\alpha \}$ are fixed,
	one would equivalently need to estimate from above the quantities 
	$\int_{U_\alpha} |\til F_\alpha|^q$ 
	and 
	$\int_{U_\alpha} \| \na_\rho  \til F_\alpha \|_\rho^q$ $\ \forall \alpha$.
	Instead of integrating $\til F_\alpha$ and $\nabla \til F_\alpha$ over $U_\alpha$, we can integrate over compact sets $V_\alpha \Subset U_\alpha$ that contain $\supp \eta_\alpha$. 
	Also, the covering $\{(U_\alpha, \vphi_\alpha)\}$ and $\{ \eta_\alpha \}$ are fixed, and $\eta_\alpha$, $\nabla \eta_\alpha$ are bounded on $V_\alpha$,
	so computing in local coordinates, it would be enough to find estimates from above of 
	$\int_{\vphi_\alpha(V_\alpha)} | F_\alpha|^q$ 
	and 
	$\int_{\vphi_\alpha(V_\alpha)} \| \na  F_\alpha \|^q$ 
	for all $\alpha$.\\
	
	

	Since there is a finite number of sets in the cover, there is such $b>0$, that for all $\alpha$, (taking suitable $\delta'$ to be the maximum of all $\delta$ taken for each $\alpha$)
	$$ \vphi_\alpha (V_\alpha) \subseteq
	\mathbf{P} =
	\{ ( z_1,\ldots, z_{2n} ) ~:~
	|z_1|,\ldots,|z_d| \leq b, ~  
	\sqrt{ z_{d+1}^2 + \ldots + z_{2n}^2 } \leq \delta'		
	\} \;.
	$$

	
	We need to check that $\int_{V_\alpha} |F_\alpha|^q$ can be made small by the same constructions. Indeed, we have
	\begin{align*}
	\int_{\vphi_\alpha(V_\alpha)} |{F}_\alpha| ^q 
	& \leq
	\int_{|z_1|,\ldots. |z_d| \leq b} dVol_{z_1,\ldots, z_d} \cdot
	\int_{\sqrt{z_{d+1}^2 + z_{2n}^2}\leq \delta'}
	|{f}_\alpha (z_{d+1},\ldots,z_{2n})| ^q ~
	dVol_{z_{d+1},\ldots,z_{2n}}
	= \\
	& = C_1 \cdot \int_{\pi_2(\mathbf{P})} |f_\alpha|^q \ d\Vol_{z_{d+1},\ldots,z_{2n}} \;.
	\end{align*}

	\noindent
	Here, the right-hand-side can be made as small as needed by \cref{lem_func_fast_decay_high_dim}, the constant $C_1$ depends only on $b$, i.e. on the fixed cover.
	Similarly, we have
	\begin{align*}
	\int_{\vphi_\alpha(V_\alpha)} \| \na F_\alpha \|^q 
	& =
	\int_{\vphi_\alpha(V_\alpha)} \| \na (f_\alpha \circ \pi_2) \|^q 
	%
	\leq \\
	&
	\leq
	C_1 \cdot \int_{\pi_2(\mathbf{P})} \| \na f_\alpha\ ( z_{d+1}(x),\ldots, z_{2n}(x) )\|^q d\Vol_{z_{d+1},\ldots, z_{2n}} \;.
	\end{align*}
	
	\noindent
	By \cref{lem_func_fast_decay_high_dim} applied to $span(e_{d+1}, \ldots, e_{2n})$, the integral on the right-hand-side can be made arbitrarily small for any $q \leq 2n-d$, given that $d \leq 2n-2$.\\

	Thus, we were able to find functions $F \in C^\infty _c (M)$ with 
	$0 \leq F \leq 1$, $F \restr_{near\ X_0} = 0$, $F \restr_{near \ X_1} = 1$ with arbitrary small 
	$\| \na_\rho F\|_q$.
	Hence (by \cref{eq: bounding_PB_using_gradient}) for any $G \in C^\infty _c (M)$ with 
	$0 \leq G \leq 1$, $G \restr_{near\ Y_0} = 0$, $G \restr_{near \ Y_1} = 1$ and for any $\epsi>0$, there exists $F$ such that 
	$(F,G)\in \calF_4'$ and $\| \pois{F,G} \|_q \leq \epsi$.
	
	We conclude that $pb_4^q (X_0,X_1,Y_0,Y_1) = 0$, $\ \forall  1\leq q\leq 2n-d$ with $d = \dim X_1 \leq 2n-2$, as required.
\end{proof}


\begin{rmk} \label{example_why_cond_on_dim_pb_4_multidim_vanishes}
	Let us note that the condition on $d = \dim X_1$ cannot be omitted.\\
	As an illustration, we explore a situation where $d = 2n - 1$, that is, when $X_1$ is a hypersurface (and $q=2n-d=1$).
	Let $(M^2, \sigma)$ be a closed symplectic surface, and let $\Pi \subset M$ be a curvilinear quadrilateral with sides $X_0, Y_0, X_1, Y_1$ listed in cyclic order.
	Pick also some closed symplectic manifold $(N^{2n-2}, \tau)$, where $n \geq 2$.
	We consider the product $M \times N$ with the symplectic form $\om = \sigma \oplus \tau$, and the quadruple $X_i' = X_i \times N, Y_i' = Y_i \times N$ for $i = 0,1$.
	We claim that $pb_4^{q=1} (X_0',X_1',Y_0',Y_1') > 0$. To prove it, we imitate the proof of \cref{lem: lem_lower_bd_L_1_pb_surface_wo_bdry}.
	
	Let $(F,G) \in \calF_4' (X_0',X_1',Y_0',Y_1')$. 
	We shall find a global estimate from below for 
	$\| \pois{F,G} \|_1 = \int_{M \times N} |\pois{F,G}| \om^n$.
	Let $U$ stand either for $\Pi$ or for $M \setm \Pi$.
	Also, we denote the endpoints of $X_1$ by
	$a = X_1 \cap Y_0$ and $b = X_1 \cap Y_1$.
	
	Using Stokes theorem once, we get
	\begin{align*}
	\int_{U \times N} |\pois{F,G}| \om^n 
	&\geq 
	\left| \int_{U \times N} \pois{F,G} \om^n \right|
	= \frac{1}{n} \left| \int_{U \times N} dF \wedge dG \wedge \om^{n-1} \right|
	= \frac{1}{n} \left| \int_{\de U \times N} FdG \om^{n-1} \right| \;.
	\end{align*}
	Using again Stokes theorem,
	\begin{align*}
	\left| \int_{\de U \times N} FdG \om^{n-1} \right| 
	&
	= \left| \int_{\de \Pi \times N} FdG \om^{n-1} \right| 
	= \left| \int_{X_1 \times N} FdG \om^{n-1} \right|
	= \left| \int_{X_1\times N} dG \om^{n-1} \right| = \\
	&
	= \left| \int_{\de(X_1 \times N)} G \om^{n-1} \right|
	= \left| \int_{b\times N} G \om^{n-1}  -  \int_{a\times N} G \om^{n-1} \right| = \\
	&
	= \int_{b \times N} \om^{n-1} = \Vol_\tau (N) \;.
	\end{align*}

	Thus, we get a positive lower bound 
	\begin{align*}
	\| \pois{F,G} \|_1 
	&=
	\int_{\Pi \times N} |\pois{F,G}| \om^{n} + \int_{(M \setm \Pi) \times N} |\pois{F,G}| \om^{n}
	\geq 
	\frac{1}{n} \cdot 2 \cdot \Vol_\tau (N)\;.
	\end{align*}
	\noindent
	Hence $pb_4^{q=1} (X_0',X_1',Y_0',Y_1') \geq \frac{2}{n} \Vol_\tau (N) > 0$.
\end{rmk}






\section{$pb_4^q$ of a curve on a surface} \label{sec: pb_4^q_of_a_curve_on_a_surface}

Let $\Sig = \Sigma_g$ be a smooth connected oriented surface of genus $g \geq 0$ without boundary.
Consider a simple closed curve $\tau$ on $\Sig$. 
Here $\tau$ is the image of an embedding $\alpha : S^1 \into \Sig$.\\
Suppose $S^1$ is divided into four closed segments $\tilde{\Delta}_1,\tilde{\Delta}_2,\tilde{\Delta}_3,\tilde{\Delta}_4$ by four distinct points in $S^1$, where the segments are listed in cyclic order. 
This induces a partition of $\tau$ into four closed segments $\Delta_i = \alpha(\tilde{\Delta}_i)$.\\

We shall consider the space $C^\infty_c (\Sig)$ of smooth compactly supported functions on $\Sig$ with the $L_q$-norm ($1\leq q \leq \infty$), and discuss $pb_4^q$ with respect to this norm. 
%

\noindent
Let us introduce 
\begin{equation}
pb_4^q(\tau):=pb_4^q(\Delta_1,\Delta_2,\Delta_3,\Delta_4) \;.
\end{equation}

\begin{claim} \label{claim: pb_4_of_curve_is_well_deifned}
	$pb_4^q(\tau)$ 
	is well-defined, i.e. it does not depend on the choice of the partition $\Delta_1,\Delta_2,\Delta_3,\Delta_4$.
\end{claim}

\begin{proof} 
	Consider two configurations of four cyclically ordered points on $\tau$, $\{a_i\}_{i=1}^{4}$ and $\{a_i'\}_{i=1}^{4}$, dividing $\tau$ into segments $\{\Delta_i\}_{i=1}^{4}$ and $\{\Delta_i'\}_{i=1}^{4}$ respectively. It is enough to show that one division can be mapped to the other by a symplectomorphism of $\Sig$.\\

	On $\tau$, take a vector field $v$, so that its flow $\{\psi_t\}$ of diffeomorphisms of $\tau$ takes $\Delta_i$ to $\Delta_i'$, i.e., 
	$\psi_1(\Delta_i)=\Delta_i'$ and $\psi_0 = id$. 

	Note that any vector field $v$ on $\tau$ can be extended to a Hamiltonian vector field on $T^*\tau$, where $\tau$ is viewed as the zero section of its cotangent bundle.  
	%
	To this end, define a Hamiltonian $H:T^*\tau \to \R$ at a point $(q,p)\in T^* \tau$ to be
	$H(q,p) = p(v(q))$, where $p\in T^*_q(\tau)$. Then $\sgrad(H)=v$ on $\tau$.
	\footnote{
		In canonical local coordinates $(p,q)$ on $T^*\tau$, $\sgrad H = (-\frac{dH}{dq},\frac{dH}{dp})$, which is $(0,v(q))$ when restricted to $\tau$. 
	}
	
	Let us indeed extend the vector field $v$ we took to a Hamiltonian vector field on $T^*\tau$, denoting the corresponding Hamiltonian by $H$ and its flow by $\{\Psi_t\}\subset \Symp(T^*\tau)$.

	By Darboux-Weinstein theorem, a neighborhood $U'$ of $\tau$ in $(T^*\tau,\om_{std})$ is symplectomorphic to a neighborhood $U$ of $\tau$ in $(\Sig,\om)$, as $\tau$ is a Lagrangian submanifold of $\Sig$. Denote this symplectomorphism by $\beta:U \to U'$.
	Multiplying $H$ by an appropriate cut-off function (that equals $1$ on $\tau$), we can guarantee $H$ to have compact support in $U'$.\\ 
	Take $\tilde{\Psi}_t = \beta^{-1} \circ \Psi_t \circ \beta$. We thus get a flow of symplectomorphisms on $\Sig$.
	Note that $\tilde{\Psi}_1 \in \Symp(\Sig)$ has compact support inside $U$, and 
	$\Psi_1(\Delta_i)= \Delta_i'$ $\forall i$.
	Hence $pb_4^q (\Delta_1,\Delta_2,\Delta_3,\Delta_4) = pb_4^q (\Delta_1',\Delta_2',\Delta_3',\Delta_4')$, 
	so $pb_4^q (\tau)$ does not depend on the division of $\tau$.
	%
	%
	%
\end{proof} 


%
We now split our investigation into two parts. First, we will examine $pb_4^q(\tau)$ for a non-separating curve $\tau$, i.e., such a curve that $S \setminus \tau$ is connected. We claim that in this case, $pb_4^q (\tau)$ vanishes. 
Then, we will examine the case of a separating curve, where the situation is different, in the sense that the result would depend on the areas of the components of $\Sig \setminus \tau$ and on $q$.

%

\begin{thm}	\label{thm: thm_pb_4_non_sep_curve_on_surface}
	If $\tau \subset \Sig$ is non-separating, then $pb_4^q (\tau) = 0$ for any $1\leq q \leq \infty$.
\end{thm}


\begin{proof}
	%
	Take two points $P',P''\in \Sig$ in a small neighborhood of $\tau$, lying on different sides of $\tau$, meaning that any curve connecting $P'$ and $P''$ that stays in a small neighborhood of $\tau$ must intersect $\tau$.
	Since $\tau$ is non-separating, there exists a simple smooth curve ${\gamma_1} \subset \Sig \setm \tau$ connecting $P'$ and $P''$. 
	Continue $\gamma_1$ by a curve $\gamma_2$ that connects the points $P'$ and $P''$, with $\gamma_2$ lying inside a small neighborhood of $\tau$, so that it does not intersect $\gamma_1$ other than at their mutual end-points, and so that $\gamma_2$ intersects $\tau$ transversally at one point $P$.
	Thus, we obtain a closed curve
	$\gamma = \gamma_1 \cup \gamma_2$ 
	that intersects $\tau$ at a unique point $P$ transversally.
	
	\begin{figure} [h!]
		\centering
		\includegraphics [scale= 0.6] {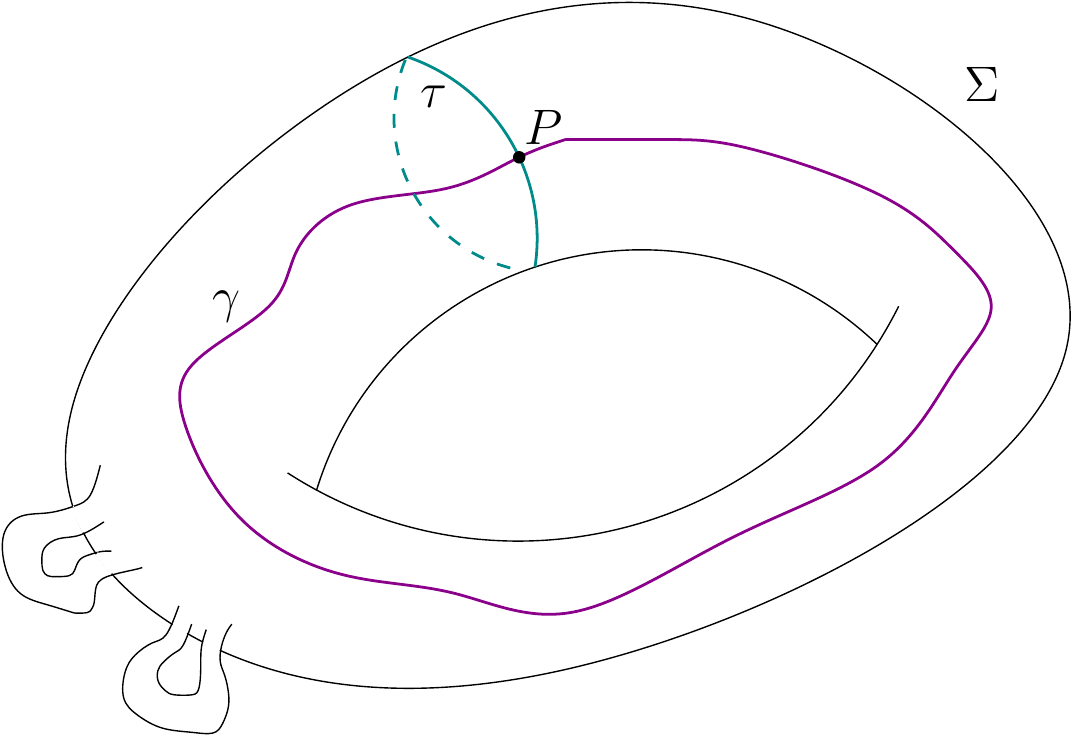}
		\; \; \; \;
		\includegraphics[scale=0.75]{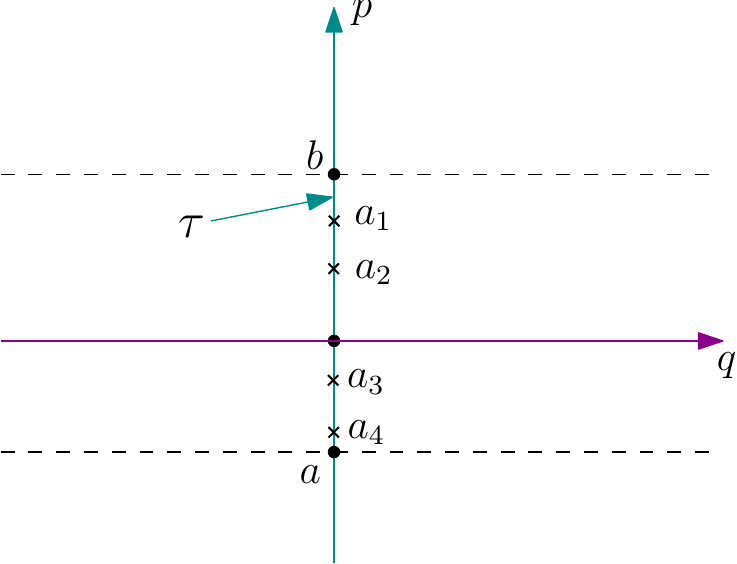}
		\caption{The curves $\tau$ and $\gamma$.}
	\end{figure}

	The curve $\gamma \subset \Sig$ is a Lagrangian submanifold, hence, by Darboux-Weinstein theorem, there exist 
	a neighborhood $U$ of $\gamma$ in $\Sig$ and
	a neighborhood $V$ of $\gamma$ in $T^* \gamma$ 
	that are symplectomorphic.
	Here we equip $T^* \gamma$ with the standard symplectic form and identify $\gamma$ with the zero section of its cotangent bundle.
	%
	
	%
	%
	%
	
	For the sake of clarity, let us indeed identify $U$ with $V$, and thus consider local coordinates $q,p$ on $U$, so that $\gamma = \{p = 0\}$. Also, without loss of generality, suppose that $U$ in these coordinates is a strip $U = \{ p\in I = (a,b) \}$, where $(a,b) \ni 0$, and $\tau \cap U = I \times \{0\}$.
	
	
	Pick four points $a < a_1 < a_2 < a_3 < a_4 < b$ on $\tau$, dividing $\tau$ into four segments
	$\Delta_i = [a_i, a_{i+1}]$ for $i = 1,2,3$ and 
	$\Delta_4$ being the closure of $\tau \setm \cup_{i=1}^3 \Delta_i$.

	Let us define a pair of functions $F,G \in \calF_4' (\Delta_1,\Delta_2,\Delta_3,\Delta_4)$. First, define them on $U$ as functions of the coordinate $p$.
	
	Consider two functions $f,g : \tau \to [0,1]$ defined as follows.
	Let $f$ be $0$ on $\Delta_1$ and $1$ on $\Delta_3$, increasing on $\Delta_2$ and decreasing on $\Delta_4$, such that it is zero outside $I = (a,b)$.
	Let $g$ be instead $0$ on $\Delta_2$ and $1$ on $\Delta_4$, increasing on $\Delta_3$ and decreasing on $\Delta_1$.
	
	Take $F(q,p) = f(p)$ and $G(q,p) = g(p)$ on $U$. 
	Further, 
	extend $F$ by zero outside $U$, and extend $G$ by $1$ outside $U$ to obtain two smooth functions defined on $\Sig$. In order for $G$ to have a compact support too, multiply it by a cutoff function that equals $1$ on $U$ and has compact support.
	Then indeed
	$(F,G) \in \calF_4' (\Delta_1,\Delta_2,\Delta_3,\Delta_4)$, and 
	$\pois{F,G} \equiv 0$, since in a neighborhood of every point on $\Sig$ either $F$ or $G$ is constant. 
	
	Hence, $pb_4^q (\tau) = 0$ for any $1 \leq q \leq \infty$.
\end{proof}

We turn now to the case when $\tau \subset \Sig$ is a separating curve. Let us first consider a concrete example as an illustration to the more general case to follow.

\begin{exm} \label{exm: exm_cylinder_A_B}
	Consider a cylinder 
	$Z_{A,B} = (0, \frac{A+B}{2\pi}) \times S^1$ with coordinates $(t,\theta)$ and the 
	area form $dt \wedge d\theta$.
	Let $\tau = \{ \frac{A}{2\pi} \} \times S^1 \subset Z_{A,B}$ be a smooth closed curve, it divides our cylinder into two components, leaving the union $Z_{A,B} \setm \tau = \Sig_1 \cup \Sig_2$ of two open subsurfaces of areas $A$ and $B$. In order to compute $pb_4^q (\tau)$ we shall map $Z_{A,B}$ symplectically to a plane region. 
	
	We denote by $B_r \subset \R^2$ an open ball of radius $r$ centered at the origin.
	Consider $M = B_{R_1} \setm B_\epsi$ for some small fixed $\epsi > 0$ and a circle
	$\tau' = \de B_{R_2}$. 
	In order to have 
	$\Area (B_{R_1} \setm B_\epsi) = B$ and
	$\Area (B_{R_2} \setm B_\epsi) = A$, 
	we put 
	$R_1^2 = \frac{B}{\pi} + \epsi^2$ and
	$R_2^2 = \frac{A}{\pi} + \epsi^2$.

	Consider the coordinates $\left( \rho = \frac{r^2}{2}, \alpha \right)$ on $M$, where $(r, \alpha)$ are polar coordinates in the plane. Equip $M$ with the area form 
	$d \rho \wedge d \alpha = r dr d\alpha$.
	
	\begin{figure} [h!]
		\centering
		\includegraphics[scale=0.75]{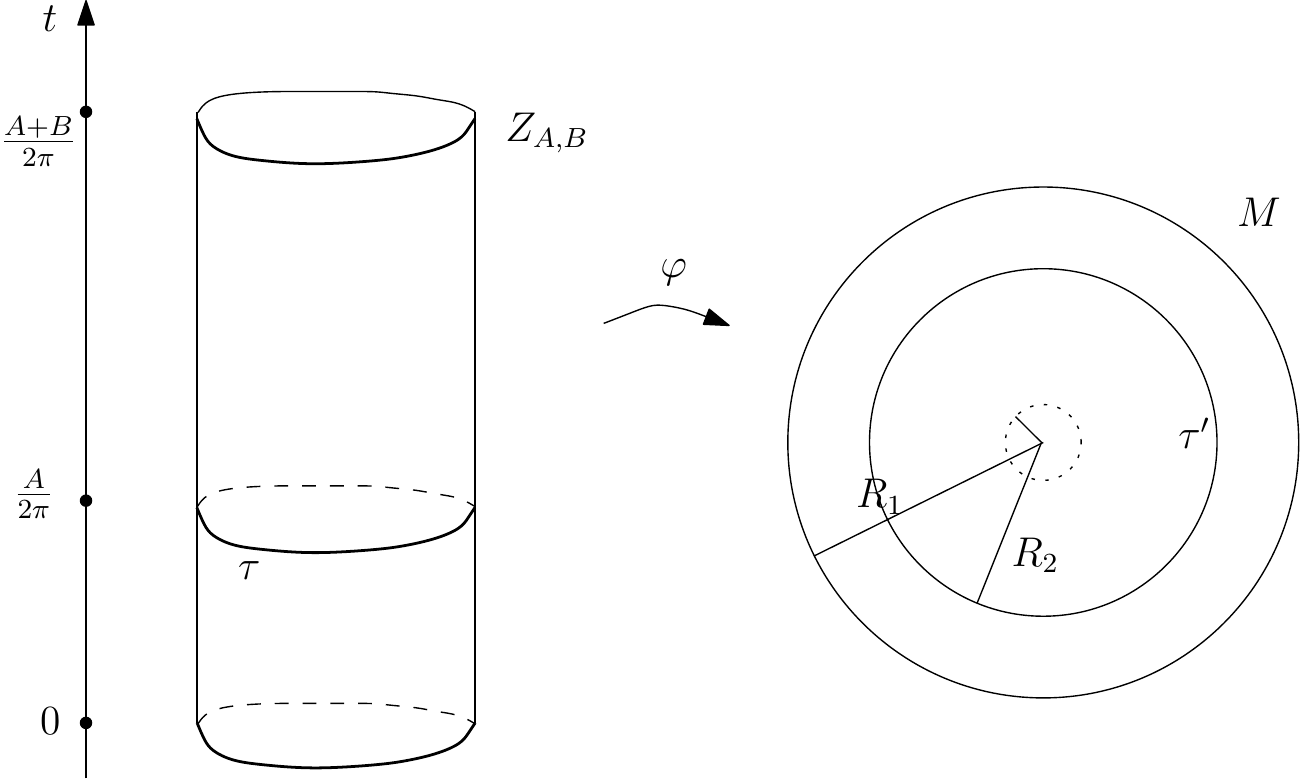}
		\caption{The cylinder $Z_{A,B}$ mapped to an annulus $M$ in the plane.}
	\end{figure}
	
	Take a map $\vphi : Z_{A,B} \to M$, defined by 
	$(t,\theta) \mapsto (\frac{\epsi^2}{2} + t, \theta)$.
	Note that $\vphi$ is indeed a symplectomorphism that takes $Z_{A,B}$ to $M$ (and $\tau$ to $\tau'$).
	Hence $pb_4^q(\tau, Z_{A,B}) = pb_4^q (\tau', M)$ for any $1 \leq q \leq \infty$, so that
	\begin{equation}
	pb_4^q(\tau, Z_{A,B}) = 	
	\begin{cases}
	2							 & \text{if } \ q=1  \\
	\left(\frac{1}{A^{q-1}} + \frac{1}{B^{q-1}} \right) ^{1/q}       & \text{if } \ 1< q < \infty \\
	\max (\frac{1}{A}, \frac{1}{B})				     & \text{if } \ q=\infty 
	\end{cases}	
	\end{equation}
	
	(similarly to what was computed in \cref{thm: thm_pb_4_q=1_dim=2}, \cref{thm: thm_q>1_dim=2_pb_4_finite_area}  , and by 
	\cite[Section 7.5.3]{polterov_rosen2014func_theory}).
\end{exm}

We can formulate the following quantitative result, claiming that $pb_4^q (\tau)$ does not vanish as long as at least one of the components of $\Sig \setminus \tau$ has finite area.
(See also \cref{rmk: separating_curve_one_component_infinite}.)

\begin{thm} \label{thm: pb_4_sep_curve_surface_no_bdry_finite_areaed_components}
	Suppose $(\Sig,\om)$ is a smooth oriented connected symplectic surface without boundary and
	$\tau \subset \Sigma$ is a separating curve. 
	Denote by $A$ the minimum of the areas of the two components of $\Sig \setm \tau$.
	If $A < \infty$, then 
	\begin{equation} \label{eq: pb_4^q_sep_curve_finite_areas_surface}
	pb_4^q(\tau) = 	
	\begin{cases}
	2							 & \text{if } \ q=1  \\
	\left(\frac{1}{A^{q-1}} + \frac{1}{B^{q-1}} \right) ^{1/q}       & \text{if } \ 1< q < \infty \\
	\max (\frac{1}{A}, \frac{1}{B})				     & \text{if } \ q=\infty 
	\end{cases}	
	\end{equation}
\end{thm}

\begin{proof}
	For any $q$, to prove that the stated values of $pb_4^q$ represent lower bounds for each case of $q$, we may use the same technique as in the proofs of lower bounds in \cref{thm: intro_thm_pb_4_dim=2_all}.
	Let us elaborate.\\
	Take some partition $\Delta_1, \Delta_2, \Delta_3, \Delta_4$ of $\tau$.
	We consider the two components $\Sigma_1, \Sigma_2$ of $\Sigma \setm \tau$.
	Take any $(F,G) \in \calF_4' (\Delta_1, \ldots, \Delta_4)$.\\
	For $q=1$, by an argument as in the proof of \cref{lem: lem_lower_bd_L_1_pb_surface_wo_bdry}, we have
	$$ 
	\int_{\Sigma_i} |\pois{F,G}| \om \geq 1, \; i=1,2 \;.
	$$
	Hence $\| \pois{F,G}\|_1 \geq 2$ for any such pair $(F,G)$, 
	so $pb_4^{q=1} (\tau) \geq 2$.
	(The same proof goes through if $\Sigma_1$ or $\Sigma_2$ (or both) have infinite area.)\\
	
	For $1 < q < \infty$, we use the case $q=1$ and H\"{o}lder inequality, imitating the proof the lower bound in \cref{thm: thm_q>1_dim=2_pb_4_finite_area} as follows.
	By the mentioned considerations, 
	we get that for any $(F,G) \in \calF_4' (\Delta_1, \ldots, \Delta_4)$,
	$$
	\int_{\Sigma_1} |\pois{F,G}|^q \om \geq \frac{1}{A^{q-1}} \;, \;
	\int_{\Sigma_2} |\pois{F,G}|^q \om \geq \frac{1}{B^{q-1}} \;.
	$$
	Hence 
	$pb_4^q (\tau) \geq \left( \frac{1}{A^{q-1}} + \frac{1}{B^{q-1}} \right)^{1/q}$.\\
	
	
	%

	
	To prove the upper bound, we will use an analogue of \cref{lem: lem_symp_stdrd_surface_w_curv_quadr}, 
	claiming that for any two numbers $0 < A' < A$ and $0 < B' < B$, 
	there is an area preserving map 
	$\vphi : Z_{A',B'} \to M$, such that it takes the circle $\sigma = \{\frac{A'}{2\pi}\} \times S^1$ to $\tau$.
	Then by \cref{lem: symp_preserves_q_norm_of_PB_any_q_g_1} applied to the symplectomorphism 
	$ \vphi : Z_{A',B'} \to \vphi (Z_{A',B'})$,
	we get that 
	\begin{equation} 
	pb_4^q(\tau) \leq 	
	\begin{cases}
	2							 & \text{if } \ q=1  \;, \\
	\left(\frac{1}{(A')^{q-1}} + \frac{1}{(B')^{q-1}} \right) ^{1/q}       & \text{if } \ 1< q < \infty \;, \\
	\max (\frac{1}{A'}, \frac{1}{B'})				     & \text{if } \ q=\infty \;. 
	\end{cases}	
	\end{equation}
	\noindent
	By an argument similar to the proof of the upper bound in \cref{thm: thm_q>1_dim=2_pb_4_finite_area}, taking $A' \to A$ and $B' \to B$, we obtain the declared result.
	
\end{proof}

\begin{rmk} \label{rmk: separating_curve_one_component_infinite}
	Using the same argument as in the proof of 
	\cref{thm: thm_pb_4_q>1_dim=2_infinite_area}, 
	we can conclude from \cref{thm: pb_4_sep_curve_surface_no_bdry_finite_areaed_components} that in case one of the areas of the components $\Sig_1, \Sig_2$ is infinite, and the other is finite (say $A<\infty$), then
	still $pb_4^q (\tau)$ is positive and 
	\begin{equation}
	pb_4^q(\tau) =
	\begin{cases}
	2							 & \text{if } \ q=1  \\
	\left(\frac{1}{A^{q-1}} \right) ^{1/q}       & \text{if } \ 1< q < \infty \\
	\frac{1}{A}				     & \text{if } \ q=\infty \;.
	\end{cases}
	\end{equation}	\\
\end{rmk}

%
%
%





\section{Discussion}

Following our results, there are some questions that require further exploration.

First, it would be interesting to complete the examination of the following functional 
(for $1 \leq p,q \leq \infty$):
\begin{align*}
\Psi_{p,q} : C^\infty_c(M) \times C^\infty_c(M) \to \R_{\geq 0},\ 
(F,G) \mapsto \liminf_{\ol F,\ol G \xrightarrow[{L_q}]{} F,G} \| \pois{\ol F,\ol G} \|_p \;.
\end{align*}
Recall that by $C^0$-rigidity of the Poisson bracket we know that 
$\Psi_{\infty,\infty} (F,G) = \| \pois{F,G} \|_\infty$, 
and by \cref{thm: thm_pb_non_rigidity_L_p_L_q_combined_intro}
$\Psi_{p,q}$ vanishes identically
for $1 \leq q<\infty$ and any $1\leq p \leq \infty$.
It would be interesting to find out whether in the remaining case $q=\infty,\ p<\infty$ this functional exhibits any rigidity.\\
To say a few words in this direction, let us recall a result by Zapolsky 
(see \cite{zapolsky_quasi_and_pb_surfaces_07})
which gives a lower bound to $\| \pois{F,G} \|_1$ in terms of the $C^0$-continuous functional 
$$
\Pi(F,G) := |\zeta(F+G) - \zeta(F) - \zeta(G)| \;
$$
that measures the non-linearity of a fixed quasi-state $\zeta$ on $M$.
The result states that for any simple quasi-state $\zeta$ on a closed symplectic surface $(M,\om)$, we have
$$
\| \pois{F,G} \|_1 \geq \Pi(F,G)^2 \;.
$$ 
Thus, we can conclude immediately that if for a pair $F,G\in C^\infty_c(M)$ we have $\Pi(F,G) > 0$, then also $\Psi_{p=1,q=\infty} (F,G) \geq \Pi (F,G)> 0$.
(Note also that using H\"older inequality, this lower bound and conclusion can be generalized to any $1 < p < \infty$.)\\
Slightly modifying the proofs in Section 3.3 of \cite{zapolsky_quasi_and_pb_surfaces_07}, one can readily show positivity of $\Psi_{p=1,\ q=\infty}$ for the case of any two-dimensional symplectic manifold, i.e. for any non-commuting $F,G\in C^\infty_c (M)$, we have 
$\liminf_{\ol F,\ol G \xrightarrow[{L_q}]{} F,G} \| \pois{\ol F,\ol G} \|_p > 0$. 

\noindent
In fact, using still the ideas in Section 3.3 of \cite{zapolsky_quasi_and_pb_surfaces_07}, one can show that in the 2-dimensional case the functional $\Psi_{p,q}$ for $1\leq p < \infty$, $q=\infty$ is lower-semicontinuous. See \cite{KS_FZ_rigidity_Lp_2016}.\\

Another question arises concerning the result about $pb_4^q (X_0, X_1, Y_0, Y_1)$ vanishing for certain quadruples in the multidimensional case (\cref{thm: intro_thm_pb_4_vanish_high_dim}). We would like to know if the condition $q \leq 2n-d$ posed on $q$ is necessary. Here the dimension $d = \dim X_1 \leq 2n-2$.\\

\section*{Acknowledgments}

This work, except for a few small additions, was written as part of the requirements for the M.Sc. degree at the School of Mathematical Sciences, Tel Aviv University.\\

\noindent
I would like to express my deepest gratitude to my advisors, Prof. Lev Buhovsky and Prof. Leonid Polterovich, for their careful guidance, immeasurable patience and their trust, 
for sharing their experience and insights, both mathematical and philosophical,
and for introducing me to the world of symplectic geometry.\\
I am indebted to Daniel Rosen for many helpful conversations and moral support,
as well as to other fellow students and friends 
for their invaluable help and encouragement. \\
Finally, I would like to thank my dear family for being there all along.\\

\bibliographystyle {plain}
\nocite{*}
\bibliography {biblio_list_2}

\begin{thebibliography}{10}

\bibitem{buhovsky_2/3_conv_pb_2010}
Lev Buhovsky.
\newblock The {$2/3$}-convergence rate for the {P}oisson bracket.
\newblock {\em Geom. Funct. Anal.}, 19(6):1620--1649, 2010.

\bibitem{buhovsky_entov_polterovich_2012poisson_brck}
Lev Buhovsky, Michael Entov, and Leonid Polterovich.
\newblock Poisson brackets and symplectic invariants.
\newblock {\em Selecta Math. (N.S.)}, 18(1):89--157, 2012.

\bibitem{cairns1961}
Stewart~S. Cairns.
\newblock A simple triangulation method for smooth manifolds.
\newblock {\em Bull. Amer. Math. Soc.}, 67:389--390, 1961.

\bibitem{cardin_viterbo_08}
Franco Cardin and Claude Viterbo.
\newblock Commuting {H}amiltonians and {H}amilton-{J}acobi multi-time
  equations.
\newblock {\em Duke Math. J.}, 144(2):235--284, 2008.

\bibitem{dacorogna_moser_1990}
Bernard Dacorogna and J{\"u}rgen Moser.
\newblock On a partial differential equation involving the jacobian
  determinant.
\newblock In {\em Annales de l'IHP Analyse non lin{\'e}aire}, volume~7, pages
  1--26, 1990.

\bibitem{entov_polterov_c0_rig_poiss_br}
Michael Entov and Leonid Polterovich.
\newblock {$C^0$}-rigidity of poisson brackets.
\newblock In {\em Symplectic topology and measure preserving dynamical
  systems}, volume 512 of {\em Contemp. Math.}, pages 25--32. Amer. Math. Soc.,
  Providence, RI, 2010.

\bibitem{mcduffSalamon_intro98}
Dusa McDuff and Dietmar Salamon.
\newblock {\em Introduction to symplectic topology}.
\newblock Oxford Mathematical Monographs. The Clarendon Press, Oxford
  University Press, New York, second edition, 1998.

\bibitem{polterov_rosen2014func_theory}
Leonid Polterovich and Daniel Rosen.
\newblock {\em Function theory on symplectic manifolds}, volume~34 of {\em CRM
  Monograph Series}.
\newblock American Mathematical Society, Providence, RI, 2014.

\bibitem{KS_FZ_rigidity_Lp_2016}
Karina Samvelyan and Frol Zapolsky.
\newblock Rigidity of the ${L}^p$-norm of the {P}oisson bracket on surfaces.
\newblock 2016.
\newblock Preprint, arXiv:1609.08891.

\bibitem{zapolsky_quasi_and_pb_surfaces_07}
Frol Zapolsky.
\newblock Quasi-states and the {P}oisson bracket on surfaces.
\newblock {\em J. Mod. Dyn.}, 1(3):465--475, 2007.

\end{thebibliography}

\bigskip \bigskip \bigskip

\noindent
\texttt{karina.samvelyan@gmail.com}


\end {document}